\documentclass{article}
\usepackage[affil-it]{authblk}

\usepackage{graphicx}
\usepackage{amsmath}
\usepackage{amsfonts}
\usepackage{amssymb}
\usepackage[amssymb]{SIunits}
\usepackage{subfigure}
\usepackage{color}
\usepackage{float}

\newcommand{\V}[1]{\ensuremath{\mathbf{#1}}}
\newcommand{\Vg}[1]{\boldsymbol{#1}}
\newcommand{\GM}[1]{\ensuremath{\mathbf{#1}}}           % Matrix

\title{On the convergence rate of the Dirichlet--Neumann iteration for unsteady thermal fluid structure interaction}

\author{Azahar Monge$^*$ and Philipp Birken%
  \thanks{e-mail: \texttt{azahar.monge@na.lu.se}; web page: http://www.maths.lu.se/staff/azahar-monge}}
\affil{$^*$Centre for Mathematical Sciences,\\ Lund University,\\ Box 118, 22100, Lund, Sweden}

\begin{document}
\maketitle

\begin{abstract}
\noindent 
We consider the Dirichlet-Neumann iteration for partitioned simulation of thermal fluid-structure interaction, also called conjugate heat transfer. We analyze its convergence rate for two coupled fully discretized 1D linear heat equations with jumps in the material coefficients across these. These are discretized using implicit Euler in time, a finite element method on one domain, a finite volume method on the other one and variable aspect ratio. We provide an exact formula for the spectral radius of the iteration matrix. This shows that for large time steps, the convergence rate is the aspect ratio times the quotient of heat conductivities and that decreasing the time step will improve the convergence rate. Numerical results confirm the analysis and show that the 1D formula is a good estimator in 2D and even for nonlinear thermal FSI applications.
\end{abstract}

\hspace*{3,6mm}\textit{Keywords:} Thermal Fluid Structure Interaction, Coupled Problems, Transmission Problem, Fixed Point Iteration, Dirichlet-Neumann Iteration % Keywords

\vspace{10pt} % Some vertical space between the abstract and first section

\section{Introduction}

The Dirichlet-Neumann iteration is a basic method in both domain decomposition and fluid structure interaction (FSI). In the latter case, the iteration arises in a partitioned approach \cite{farhat:04}, where different codes for the sub-problems are reused and the coupling is done by a master program which calls interface functions of the other codes. This allows to reuse existing software for each sub-problem, in contrast to a monolithic approach, where a new code is tailored for the coupled equations. To satisfy coupling conditions at the interface, the subsolvers are iterated by providing Dirichlet- and Neumann data for the other solver in a sequential manner, giving rise to its name.

In the domain decomposition context, the iteration has two main problems, namely slow convergence and the need for an implementation using a red-black colouring. The slow convergence can be slightly improved using a relaxation procedure. In fluid structure interaction, there are typically only two domains, coupled along an interface, making the application straight forward. The convergence rate for the interaction of a flexible structure with a fluid has been analyzed in \cite{vanBrummelen:09}. There, the added mass effect is proven to be dependent on the step size for compressible flows and independent for incompressible flows. However, the convergence rate is not great for the coupling between a compressible fluid and a structure \cite{deparis:03}, which is why a lot of effort goes into convergence acceleration. Furthermore, for incompressible fluids it is known that the ratio of densities of the materials plays an important role \cite{badia:08,causin:05}. Finally, the Dirichlet-Neumann iteration was reported to be a very fast solver for thermal fluid structure interaction \cite{biglkm:15}. 

Our prime motivation here is thermal interaction between fluids and structures, also called conjugate heat transfer. There are two domains with jumps in the material coefficients across the connecting interface. Conjugate heat transfer plays an important role in many applications and its simulation has proved essential \cite{banka:05}. Examples for thermal fluid structure interaction are cooling of gas-turbine blades, thermal anti-icing systems of airplanes \cite{buchli:10}, supersonic reentry of vehicles from space \cite{mehta:05,hinrad:06}, gas quenching, which is an industrial heat treatment of metal workpieces \cite{hefiba:01,stshle:06} or the cooling of rocket nozzles \cite{kohoha:13,kotirh:13}.

For the case of coupled heat equations, a 1D stability analysis was presented by Giles \cite{giles:97}. There, an explicit time integration method was chosen with respect to the interface unknows. On the other hand, Henshaw and Chand provided in \cite{Henshaw:2009} a method to analyze stability and convergence speed of the Dirichlet-Neumann iteration in 2D based on applying the continuous Fourier transform to the semi-discretized equations. They show that the ratios of thermal conductivities and diffusivities of the materials play an important role. This is similar to the result names above in classical FSI with incompressible fluids where the performance is affected by the ratio of densities of the materials \cite{badia:08,causin:05}. 

However, in the fully discrete case we observe in some cases that the iteration behaves differently, because some aspects of the problem are not taken into account by the semidiscrete analysis: The effect of $\Delta t$ is not accurately represented and neither are possibly different mesh widths in the two problems. This matters particularly for compressible fluids where a high aspect ratio grid is needed to accurately represent the boundary layer. This leads to geometric stiffness that significantly influences the convergence rate, as we will show. 

For the fully discrete case, the convergence rate is in principle analyzed in any standard book on domain decomposition methods, e.g. \cite{quaval:99, Toselli:2004qr}. There, the iteration matrix is derived in terms of the stiffness and mass matrices of finite element discretizations and the convergence rate is the spectral radius of that. However, this does not provide a quantitative answer, since the spectral radius is unknown. Computing the spectral radius is in general a non trivial task. In our context, the material properties are discontinuous across the interface and as a consequence, computing the spectral radius of the iteration matrix is even more difficult. 

In \cite{Mongelic:16,Monge:16}, a convergence analysis of the Dirichlet-Neumann iteration for the unsteady transmission problem using finite element methods (FEM) on both subdomains is presented. A similar analysis using finite differences (FDM) on one domain and FEM on the other one can be found in \cite{Monge:2016}. In addition, the corresponding analysis when coupling finite volumes (FVM) with FEM is described in \cite{Birken:2016,Mongelic:16}. All these results assume equal mesh sizes on both subdomains, i.e, the aspect ratio is equal to one.  

Thus, we present here a complete discretization of the coupled problem using FVM in space on one domain and FEM on the other one with variable aspect ratio $r$. We consider this to be a relevant case, because these are the standard discretizations for the subproblems. The implicit Euler method is used for the temporal discretization. Then, we derive the spectral radius of the iteration matrix exactly in terms of the eigendecomposition of the resulting matrices for the one dimensional case. The asymptotic convergence rates when approaching the continuous case in either time or space are also determined. In the spatial limit, the convergence rate turns out to be proportional to the aspect ratio $r$, whereas in the temporal limit, we obtain 0. Note that for FEM-FEM couplings, this is not the case. Moreover, we also include numerical results where it is shown that the one dimensional formula is a good estimator for a 2D version of the coupled heat equations and for two non linear FSI models, namely the cooling of a flat plate and the cooling of a flanged shaft.     

An outline of the paper now follows. In section 2, we describe the model and discretization, as well as the coupling conditions and the Dirichlet-Neumann iteration. Two thermal FSI test cases are introduced in section 3: the cooling of a flat plate and of a flanged shaft. For these, we present numerical convergence rates, motivating further analysis. A model problem, consisting of two coupled discretized heat equations, is presented in section 4 and then analyzed in 1D in section 5. In section 6, extension of the analysis to 2D and different discretizations are discussed. In section 7, the analytical results are compared to linear and nonlinear numerical results. 

%---------------------------------------------------------------------------------------

\section{Thermal FSI Methodology}

The basic setting we are in is that on a domain $\Omega_1 \subset
\mathbb{R}^d$ where $d$ corresponds to the spatial dimension, the physics is described by a fluid model, whereas on a
domain $\Omega_2 \subset \mathbb{R}^d$, a different model describing the structure is used. The two domains are almost disjoint in that they
are connected via an interface. The part of the interface where the
fluid and the structure are supposed to interact is called the
coupling interface 
$\Gamma \subset \partial \Omega_1 \cup \partial \Omega_2$. Note that $\Gamma$ might be a true subset of the intersection, because the structure could be insulated. At the interface $\Gamma$, coupling conditions are prescribed that model the interaction between fluid and structure. For the thermal coupling problem, these conditions are that temperature and the normal component of the heat flux are continuous across the interface. 

\subsection{Fluid Model}
We model the fluid using the time dependent compressible Navier-Stokes
equations, which are a second order system of conservation laws
(mass, momentum, energy) modeling compressible flow. We
consider the two dimensional case, written in conservative variables
density $\rho$, momentum ${\bf m}=\rho {\bf v}$ and energy per unit
volume $\rho E$ as:
\begin{align}
\partial_t \rho + \nabla \cdot \rho{\bf v} &=0, \nonumber \\
\partial_t \rho v_i + \sum_{j=1}^2\partial_{x_j}(\rho v_i v_j + p\delta_{ij})&= \frac{1}{Re}\sum_{j=1}^2\partial_{x_j} S_{ij}, \quad i=1,2,\\
\partial_t \rho E + \nabla \cdot (\rho H v_j)&= \frac{1}{Re}\sum_{j=1}^2 \partial_{x_j} \left( S_{ij}v_i  + \frac{q_j}{Pr} \right ). \nonumber
\end{align}
Here, enthalpy is given by $H = E + p/\rho$ with $p = (\gamma -1)\rho(E-1/2|v|^2)$ being the pressure and $\gamma=1.4$ the adiabatic index for an ideal gas. Furthermore, ${\bf q}_f = (q_1,q_2)^T$ represents the heat flux and $\textbf{S} = (S_{ij})_{i,j=1,2}$ the viscous shear stress tensor. As the equations are dimensionless, the Reynolds number $Re$ and the Prandtl number $Pr$ appear. The system is closed by the equation of state for the pressure $p = (\gamma -1) \rho e$, the Sutherland law representing the correlation between temperature and viscosity, as well as the Stokes hypothesis. Additionally, we prescribe appropriate boundary conditions at the boundary of $\Omega_1$ except for $\Gamma$, where we have the coupling conditions. In the Dirichlet-Neumann coupling, a temperature value is enforced at $\Gamma$. 

%are modelled using the Spallart-Allmaras model \cite{spaall:92}. 

\subsection{Structure Model}
Regarding the structure model, we will consider heat conduction
only. Thus, we have the nonlinear
heat equation for the structure temperature $\Theta$
\vskip-.6cm
\begin{eqnarray}
  \rho({\bf x}) c_p(\Theta) \frac{d}{dt}\Theta(\V{x},t) = - \nabla \cdot \V{q}(\V{x},t),
  \label{eq:heatcond}
\end{eqnarray}
where 
\begin{equation*}
  \V{q}_s(\V{x},t) = - \lambda(\Theta) \nabla \Theta(\V{x},t)
\end{equation*}
denotes the heat flux vector. For alloys, the specific heat capacity $c_p$ and heat conductivity $\lambda$ are temperature-dependent and highly nonlinear.  

As an example, an empirical model for the steel 51CrV4 was suggested in \cite{quhrss:11}. This was obtained measurements and a least squares fit to a chosen curve. The coefficient functions are then
\vskip-.6cm
\begin{eqnarray}
\label{heatconductivitysteel}
\lambda(\Theta)=40.1 + 0.05 \Theta - 0.0001\Theta^2 + 4.9 \cdot 10^{-8}\Theta^3
\end{eqnarray}
and
\vskip-.6cm
\begin{eqnarray}
\label{heatcapacitysteel}
c_p(\Theta)=-10\ln \left( \frac{e^{-c_{p1}(\Theta)/10}+e^{-c_{p2}(\Theta)/10}}{2}\right)
\end{eqnarray} 
with 
\vskip-.6cm
\begin{eqnarray}
c_{p1}(\Theta)=34.2e^{0.0026\Theta} + 421.15
\end{eqnarray} 
and 
\vskip-.6cm
\begin{eqnarray}
c_{p2}(\Theta)=956.5e^{-0.012(\Theta-900)} + 0.45\Theta.
\end{eqnarray} 
For the mass density one has $\rho =
\unit{7836}{\kilo\gram\per\cubic\meter}$. 

Finally, on the boundary, we have Neumann
conditions $\V{q}_s(\V{x},t) \cdot \V{n}(\V{x}) =
q_b(\V{x},t)$. 

\subsection{Coupling Conditions}

As mentioned before, the coupling conditions are that temperature and the normal component of the heat flux are continuous across the interface, i.e;

\begin{align}
T (\mathbf{x},t) = \mathbf{\Theta} (\mathbf{x},t), \  \mathbf{x} \in \Gamma,
\end{align}
where $T$ is the fluid temperature and $\mathbf{\Theta}$ the structure temperature and

\begin{align}
\mathbf{q}_f(\mathbf{x},t) \cdot \mathbf{n} (x) = \mathbf{q}_s(\mathbf{x},t) \cdot \mathbf{n} (x), \  \mathbf{x} \in \Gamma.
\end{align}

\subsection{Discretization in Space}
Following the partitioned coupling approach, we discretize the two models separately in space. For the fluid, we use a finite volume method, leading to the following equation for all unknowns on $\Omega_1$:

\begin{eqnarray}
\label{eq:ODEfluid}
\frac{d}{dt} {\bf u} + {\bf h}({\bf u},\Vg{\Theta}_{\Gamma}) = \V{0},
\end{eqnarray}
where ${\bf h}({\bf u},\Vg{\Theta}_{\Gamma})$ represents the spatial
discretization and its dependence on the temperatures on the discrete interface
to the
structure, here denoted by $\Vg{\Theta}_{\Gamma}$. 

Regarding structural mechanics, the use of finite element methods is ubiquitious. Therefore, we will also follow that approach here, using quadratic finite element one gets the following nonlinear equation for all unknowns on $\Omega_2$:

\begin{eqnarray}
\label{eq:ODEheat}
\GM{M}({\bf \Theta}) \frac{d}{dt}{\bf \Theta} + \GM{A}({\bf \Theta}) {\bf \Theta}
  = {\bf q}_b^f+{\bf q}_b^{\Gamma}({\bf u}).
\end{eqnarray}
Here, $\GM{M}$ is the mass matrix, also called heat capacity matrix for this problem and $\GM{A}$ is the heat conductivity and stiffness
matrix. The vector ${\bf \Theta}$ consists of all discrete temperature
unknowns and ${\bf q}_b^{\Gamma}({\bf {\bf u}})$ is the discrete heat flux vector on the
coupling interface to the fluid, whereas $ {\bf q}_b^f$ corresponds to
boundary heat fluxes independent of the fluid, for example at insulated
boundaries. 

\subsection{Time Discretization}

In time, we use the implicit Euler method with constant time step $\Delta t$. For the system \eqref{eq:ODEfluid}-\eqref{eq:ODEheat} we obtain

\begin{align}\label{discrete-fluid}
\mathbf{u}^{n+1}-\mathbf{u}^n + \Delta t\mathbf{h} (\mathbf{u}^{n+1}, \mathbf{\Theta}_{\Gamma}^{n+1}) = \mathbf{0},
\end{align}

\begin{align}\label{discrete-solid}
\mathbf{M} (\mathbf{\Theta}^{n+1}) (\mathbf{\Theta}^{n+1} -\mathbf{\Theta}^n)+ \Delta t\mathbf{A} (\mathbf{\Theta}^{n+1}) \mathbf{\Theta}^{n+1} = \Delta t(\mathbf{q}_b^f + \mathbf{q}_b^{\Gamma} (\mathbf{u}^{n+1})).
\end{align}

\subsection{The Dirichlet-Neumann Method}

The Dirichlet-Neumann method is a basic iterative substructuring method in domain decomposition and it is a common choice for treating FSI problems. Therefore, we now employ it to solve the system \eqref{discrete-fluid}-\eqref{discrete-solid}. This corresponds to alternately solving equation \eqref{discrete-fluid} on $\Omega_1$ with Dirichlet data on $\Gamma$ and \eqref{discrete-solid} on $\Omega_2$ with Neumann data on $\Gamma$.

Thus, one gets for the $k$-th iteration the two decoupled equation systems

\begin{align}
\label{DNmethodeq1}
\mathbf{u}^{n+1,k+1}-\mathbf{u}^n + \Delta t\mathbf{h} (\mathbf{u}^{n+1,k+1}, \mathbf{\Theta}_{\Gamma}^{n+1,k}) = \mathbf{0},
\end{align}

\begin{align}
\label{DNmethodeq2}
\mathbf{M} (\mathbf{\Theta}^{n+1,k+1}) (\mathbf{\Theta}^{n+1,k+1} -\mathbf{\Theta}^n)+ \Delta t\mathbf{A} (\mathbf{\Theta}^{n+1,k+1}) \mathbf{\Theta}^{n+1,k+1} = \Delta t(\mathbf{q}_b^f + \mathbf{q}_b^{\Gamma} (\mathbf{u}^{n+1,k+1})),
\end{align}
with some initial condition $\mathbf{\Theta}_{\Gamma}^0$. The iteration is terminated according to the standard criterion 
\begin{equation}\label{reltolnolinear}
\| \mathbf{\Theta}_{\Gamma}^{k+1} - \mathbf{\Theta}_{\Gamma}^k \| \leq \tau
\end{equation}
where $\tau$ is a user defined tolerance.

\section{Thermal FSI Test Cases}

In this section we present two thermal FSI test cases that are solved using the methodology explained in the previous section. The aim of this paper is to estimate the convergence rates of the Dirichlet-Neumann iteration used as a solver for thermal FSI problems. Therefore, we first want to illustrate the behavior for two examples before proceeding to the convergence analysis in the next section. Two different test cases are discussed: the cooling of a flat plate and the cooling of a flanged shaft. For the first problem, structured grids are used and for the second, unstructured grids.    

For the coupling, the Dirichlet-Neumann method as presented in \eqref{DNmethodeq1}-\eqref{DNmethodeq2} is used. A fixed tolerance of $1e-8$ is chosen for all involved equation solvers. The coupling code used has been developed in a series of papers \cite{biquhm:11,biquhm:10,biglkm:15}. It's main feature is time adaptivity, which is not employed here. The coupling between the solvers is done using the Component Template Library (CTL) \cite{manist:06}. In the fluid, the DLR TAU-Code in its 2014.2 version is employed \cite{gerfeg:97}, which is a cell-vertex-type finite volume method with AUSMDV as flux function and a linear reconstruction to increase the order of accuracy. The finite element code uses quadratic finite elements and is the inhouse code Native of the Institute for Static and Dynamic at the University of Kassel. 

\subsection{Flow over a plate}

\begin{figure}[h]
	\centering
	\includegraphics[width=0.75\textwidth]{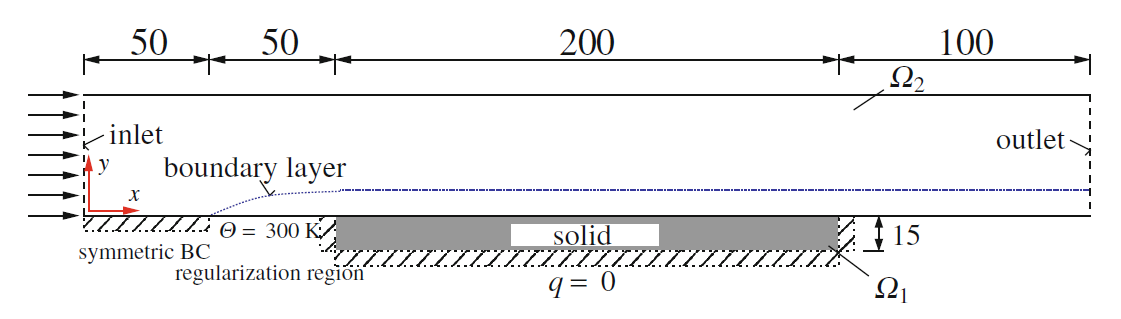}
	\caption{Sketch of the cooling of a flat plate.}
	\label{fig:hotplate}
\end{figure}

The first test case is the cooling of a flat steel plate resembling a simple work piece \cite{biquhm:11}. The work piece is initially at a much higher temperature than the fluid and then cooled by a constant laminar air stream, see figure \ref{fig:hotplate}.

The inlet is given on left, where air enters the domain with an initial velocity of $\mbox{Ma}_{\infty} = 0.8$ in horizontal direction and a temperature of 273K. Regarding the initial condition in the structure, a constant temperature of 900K at $t=0$ is chosen throughout. 

\begin{figure}[h]
	\centering
	\includegraphics[width=0.75\textwidth]{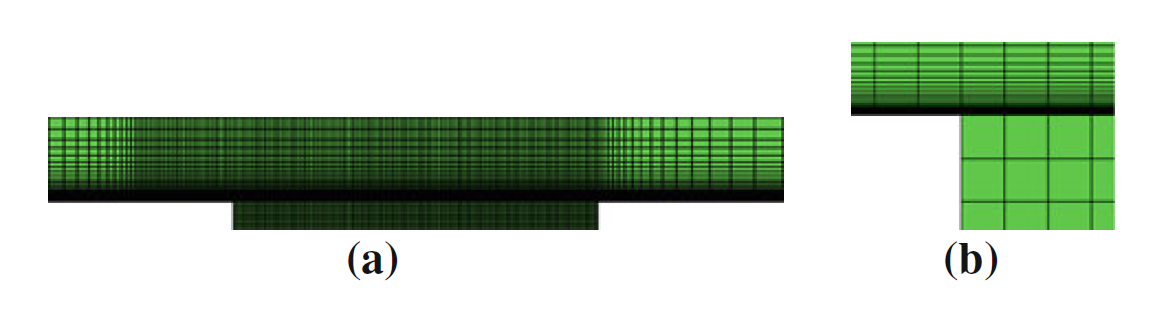}
	\caption{Full grid (left) and zoom into coupling region (right).}
	\label{fig:hotplategrid}
\end{figure}

The grid, see figure \ref{fig:hotplategrid}, is chosen cartesian and equidistant in the structural part. In the fluid region the thinest cells are on the boundary and then become coarser in $y$-direction with a maximal aspect ratio of $r=1.7780e5$. The points of the primary fluid grid and the nodes of the structural grid match on the interface $\Gamma$ and there are 9660 cells in the fluid region and $n_x \times n_y = 120 \times 9 = 1080$ elements with $121 \times 10 = 1210$ nodes in the region of the structure. 

Figure \ref{fig:coolingsystemsrates}a shows the convergence behaviour of the Dirichlet-Neumann iteration against the time step $\Delta t$. One observes how the convergence rates is roughly proportional to the time step $\Delta t$. Furthermore, even for $\Delta t=1$ a reduction of the error by a factor of ten per iteration is achieved. 

\subsection{Cooling of a flanged shaft}

\begin{figure}[h]
	\centering
	\includegraphics[width=0.5\textwidth]{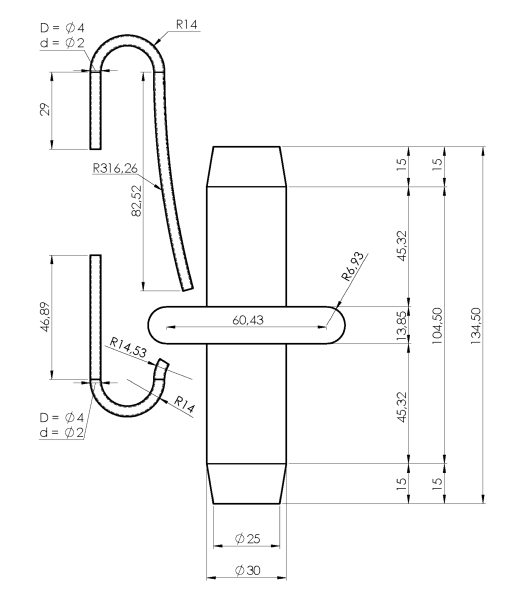}
	\caption{Sketch of the cooling of the flanged shaft.}
	\label{fig:flangedshaft}
\end{figure}

The second test case is the cooling of a flanged steel shaft by cold high pressured air (this process is also known as gas quenching) \cite{wesast:07}. Here, we have a hot flanged shaft that is cooled by cold high pressured air coming out of small tubes, see figure \ref{fig:flangedshaft}. We assume symmetry along the horizontal axis in order to consider one half of the flanged shaft and two tubes blowing air at it. We also assume that the air leaves the tube in straight and uniform way at a Mach number of 1.2. Moreover, we assume a freestream in $x$-direction of Mach 0.005. The Reynolds number is $Re = 2500$ and the Prandtl number $Pr = 0.72$.

\begin{figure}[h]
	\centering
	\includegraphics[width=0.9\textwidth]{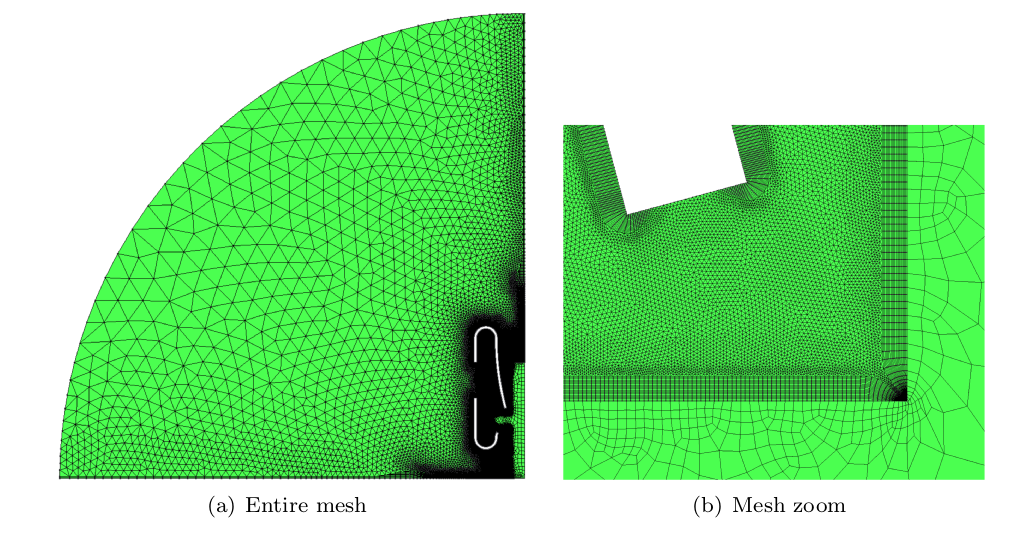}
	\caption{Full grid (left) and zoom into shaft region (right).}
	\label{fig:flangedshaftgrid}
\end{figure}

The grid, see figure \ref{fig:flangedshaftgrid}, consists of 279212 cells in the fluid, which is the dual grid of an unstructured grid of quadrilaterals in the boundary layer and triangles in the rest of the domain, and 1997 quadrilateral elements in the structure. Regarding the initial conditions, we use the procedure explained in \cite{biglkm:15}.

\begin{figure}[h!]
	\centering
	\subfigure[Test case 1: Flow over a plate]{\includegraphics[width=6cm]{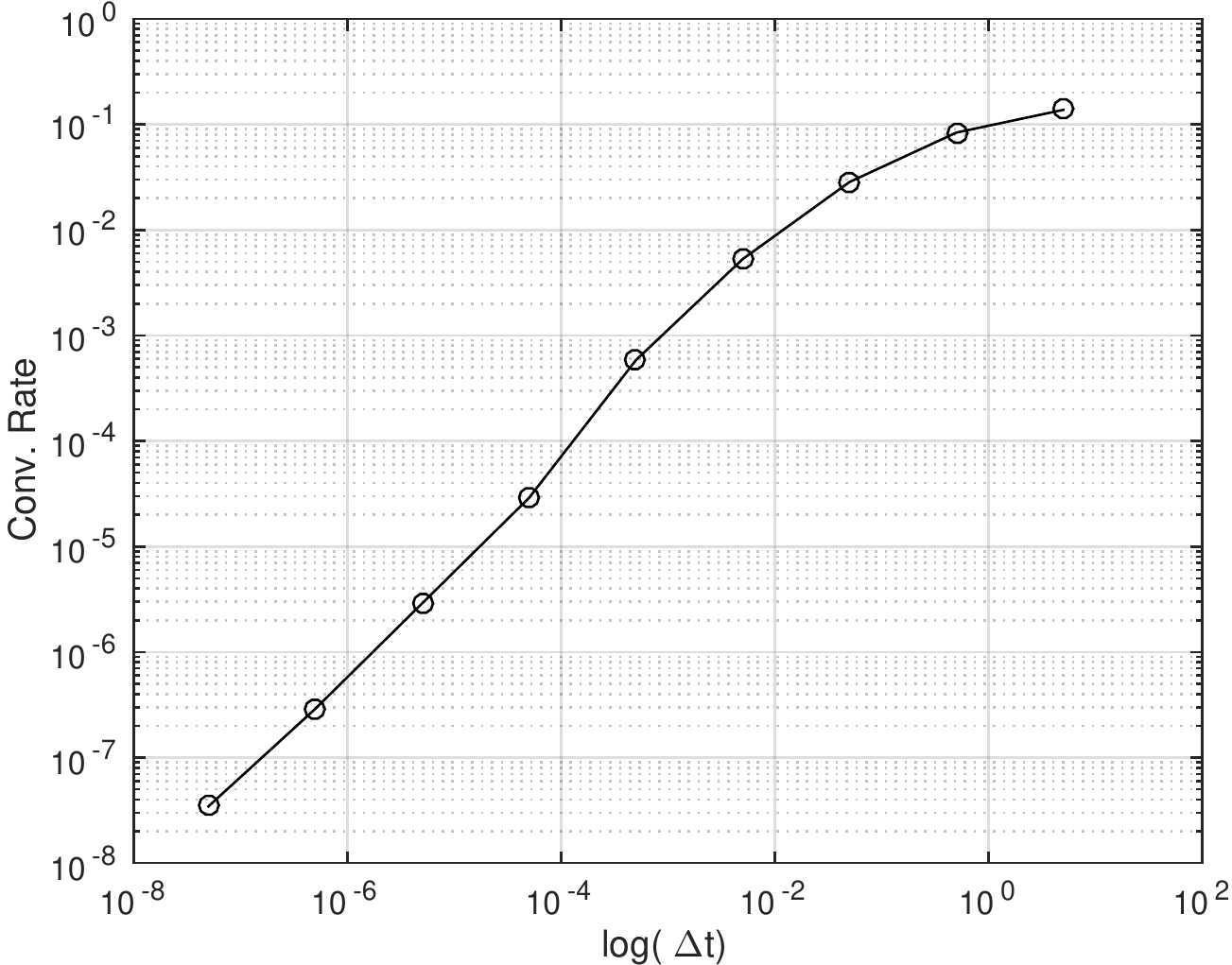}} \hfill
	\subfigure[Test case 2: Cooling of a flanged shaft]{\includegraphics[width=6cm]{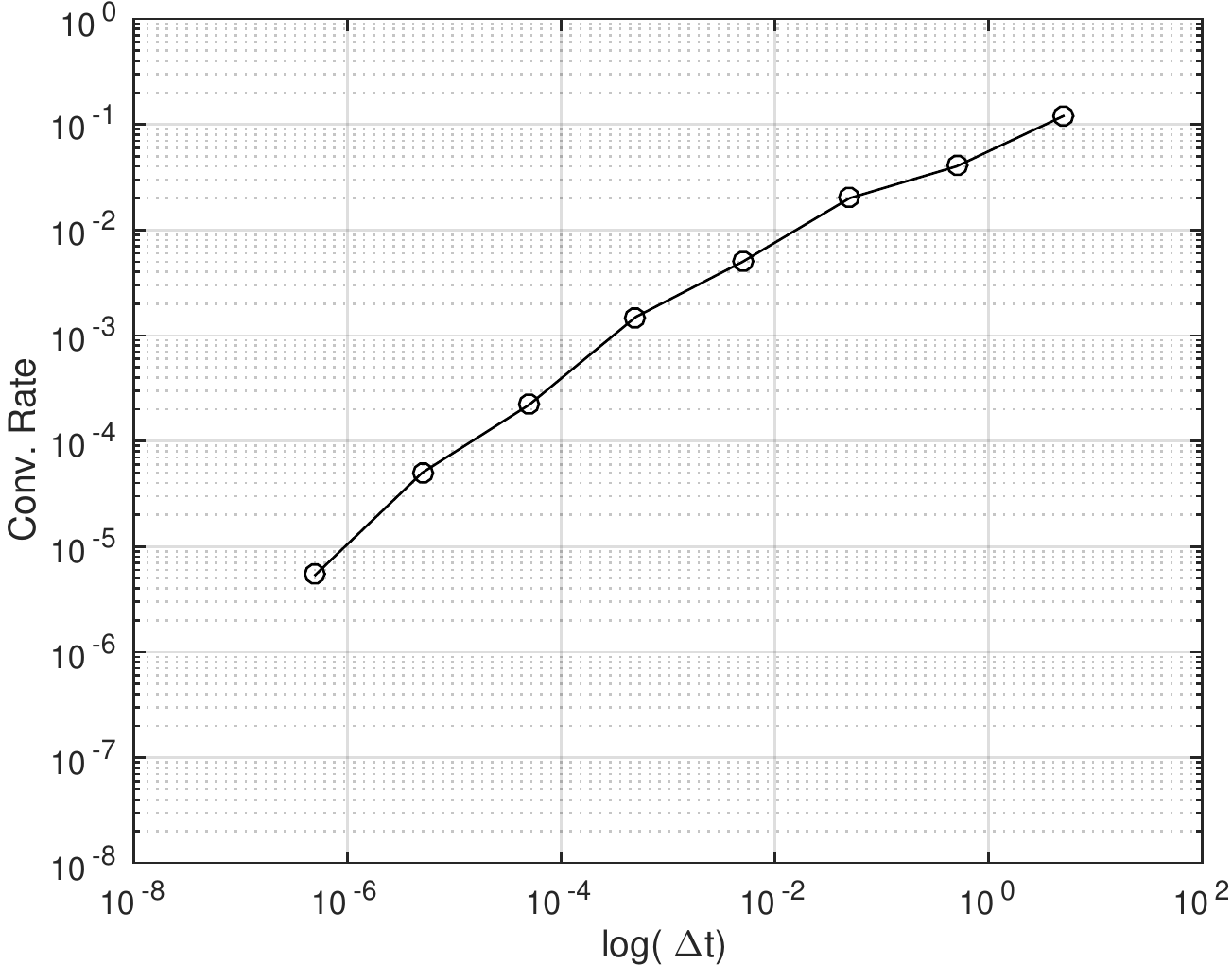}}
	\caption{Convergence behavior of the cooling systems with respect to $\Delta t$.}
	\label{fig:coolingsystemsrates}
\end{figure}

Figure \ref{fig:coolingsystemsrates}b shows the convergence behaviour of the Dirichlet-Neumann iteration against the time step $\Delta t$. The convergence rate is again about proportional to the time step size and again convergent even for very large time steps. If we compare the rates for the two problems, we observe that for a given $\Delta t$, the iteration is about a factor ten faster for the plate.

Summarizing, the Dirichlet-Neumann iteration is a very fast solver for thermal FSI. To understand this better, we perform in the next section a convergence analysis for the case of two coupled linear heat equations.

\section{A Model Problem: Coupled Heat Equations}

We present here a convergence analysis of the unsteady transmission problem with mixed discretizations. In particular, we choose a finite volume method (FVM) on the first subdomain and a finite element method (FEM) on the second subdomain.  

\subsection{Model Problem}

The unsteady transmission problem is as follows, where we consider a domain $\Omega \subset \mathbb{R}^d$ which is cut into two subdomains $ \Omega_1 \cup \Omega_2 = \Omega$ with transmission conditions at the interface $\Gamma = \Omega_1 \cap \Omega_2$:

\begin{align}
\label{model:eq}
\begin{split}
\alpha_m \frac{\partial u_m(\textbf{x},t)}{\partial t} - \nabla \cdot (\lambda_m  \nabla u_m(\textbf{x},t)) & = 0,\  \  t \in [t_0, t_f],  \  \  \textbf{x} \in \Omega_m \subset \mathbb{R}^d, \  m=1,2, \\
u_m(\textbf{x},t) & = 0, \  \  t \in [t_0, t_f], \  \  \textbf{x} \in \partial \Omega_m \backslash \Gamma, \\
u_1(\textbf{x},t) & = u_2(\textbf{x},t), \  \  \textbf{x} \in \Gamma, \\
\lambda_2 \frac{\partial u_2(\textbf{x},t)}{\partial \textbf{n}_2} & = -\lambda_1 \frac{\partial u_1(\textbf{x},t)}{\partial \textbf{n}_1}, \  \  \textbf{x} \in \Gamma, \\
u_m(\textbf{x},0) & = u_m^0(\textbf{x}), \  \  \textbf{x} \in \Omega_m,
\end{split}
\end{align}
where $\textbf{n}_m$ is the outward normal to $\Omega_m$ for $m=1,2$.

The constants $\lambda_1$ and $\lambda_2$ describe the thermal conductivities of the materials on $\Omega_1$ and $\Omega_2$ respectively. $D_1$ and $D_2$ represent the thermal diffusivities of the materials and they are defined by 

\begin{align}
D_m = \frac{\lambda_m}{\alpha_m}, \  \  \mbox{with} \  \  \alpha_m = \rho_m c_{p_m} 
\end{align}
where $\rho_m$ represents the density and $c_{p_m}$ the specific heat capacity of the material placed in $\Omega_m$, $m=1,2$. 

%In the one-dimensional case ($d=1$), we discretize this problem with a constant mesh width of $\Delta x = 1/(N_m+1)$ with $N_m$ being the number of interior space discretization points in the intervals $\Omega_m$, $m=1,2$. If instead we consider \eqref{model:eq} with $d=2$, we will use a non equidistant grid in $\Omega_1$ i.e, $\Delta x_1 \neq \Delta y_1$. On the other hand, an equidistant grid will be used in $\Omega_2$ i.e, $\Delta y_2 = \Delta x_2$. For simplicity, we assume in this case that $\Delta y := \Delta y_1 = \Delta y_2 = \Delta x_2 = 1/(N_2 +1)$ and $\Delta x := \Delta x_1 = 1/(N_1+1)$ with $N_1 \neq N_2$ as represented in \ref{fig:triangles}. 

We always use the implicit Euler method for time discretization. With regards to the spatial discretization, we use FVM on $\Omega_1$ and FEM on $\Omega_2$.

\subsection{Semidiscrete Analysis}
\label{semidiscrete-section}

Before we present in the next section an analysis for the fully discrete equations, we want to describe previous results about the behaviour of the Dirichlet-Neumann iteration for the transmission problem in the semi discrete case.

Henshaw and Chand applied in \cite{Henshaw:2009} the implicit Euler method for the time discretization on both equations in \eqref{model:eq} but kept the space continuous. Then, they applied the Fourier transform in space in order to transform the second order derivatives into algebraic expressions. Once they have a coupled system of algebraic equations, they insert one into the other and obtain the Dirichlet-Neumann convergence rate $\beta$:

\begin{align}
\label{semidiscretelimit}
\beta = \left| -\frac{\lambda_1}{\lambda_2} \sqrt{\frac{D_2}{D_1}} \frac{\tanh \left( -\frac{1}{\sqrt{D_2 \Delta t}} \right)}{\tanh \left( \frac{1}{\sqrt{D_1 \Delta t}} \right)} \right|.
\end{align} 

For $\Delta t$ small enough, we have $\tanh \left( -1 / \sqrt{D_2 \Delta t} \right) \approx -1 $ and $\tanh \left( 1 / \sqrt{D_1 \Delta t} \right) \approx 1 $ and therefore:

\begin{align}
\beta \approx \frac{\lambda_1}{\lambda_2} \sqrt{\frac{D_2}{D_1}}.
\end{align}
On the other hand, for $\Delta t$ big enough, we have $\tanh \left( -1 / \sqrt{D_2 \Delta t} \right) \approx -1 / \sqrt{D_2 \Delta t}$ and $\tanh \left( 1 / \sqrt{D_1 \Delta t} \right) \approx 1 / \sqrt{D_1 \Delta t}$ and therefore:

\begin{align}
\label{semidiscreteapprox}
\beta \approx \frac{\lambda_1}{\lambda_2} \sqrt{\frac{D_2}{D_1}} \frac{\sqrt{D_1 \Delta t}}{\sqrt{D_2 \Delta t}} = \frac{\lambda_1}{\lambda_2}.
\end{align}

\begin{figure}[h]
	\centering
	\includegraphics[width=0.55\textwidth]{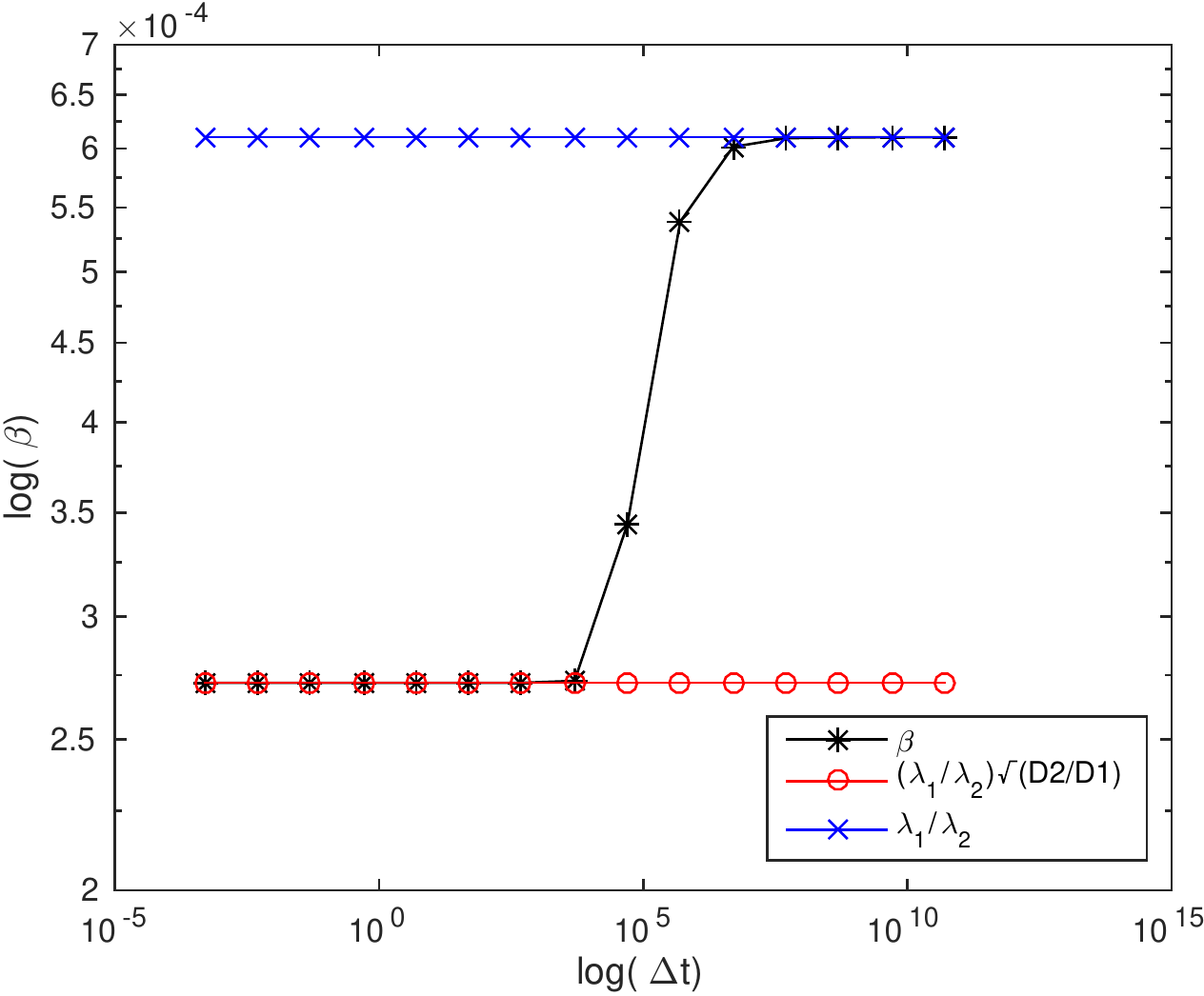}
	\caption{Semidiscrete estimator $\beta$ in \eqref{semidiscretelimit} against $\Delta t$.}
	\label{fig:semidiscretevsdeltat}
\end{figure}

Figure \ref{fig:semidiscretevsdeltat} shows $\beta$ as a function of $\Delta t$. It is almost constant, except for a short dynamic transition between $(\lambda_1/ \lambda_2) \sqrt{D_2/D_1}$ and $\lambda_1 / \lambda_2$. 

Finally, one observes in \eqref{semidiscreteapprox} that the convergence rates of the Dirichlet-Neumann iteration are given by the quotient of thermal conductivities for $\Delta t$ large. This suggests that strong jumps in the thermal conductivities of the materials give fast convergence.

\subsection{Space Discretization}

\begin{figure}[h]
	\centering
	\includegraphics[width=0.55\textwidth]{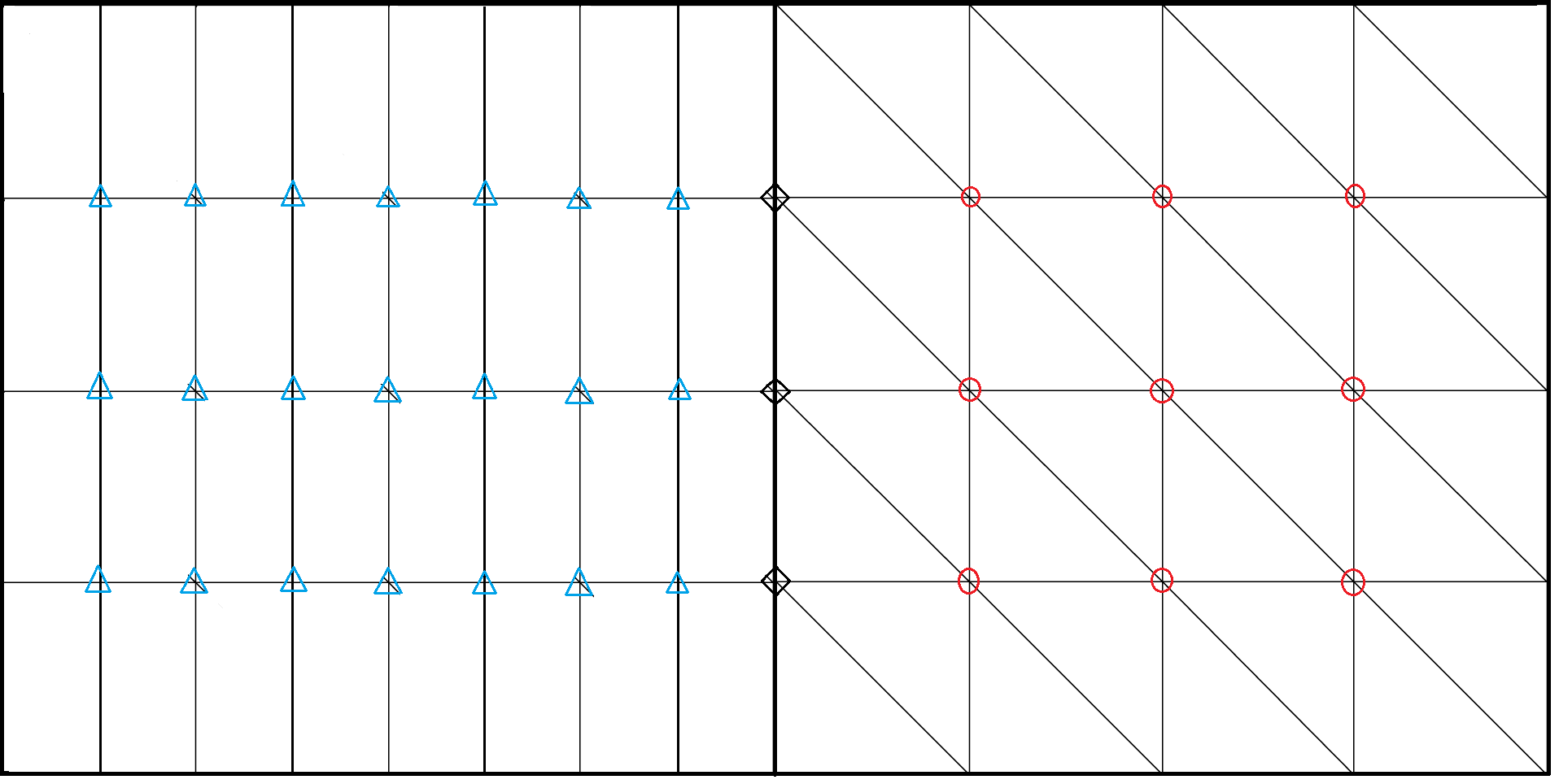}
	\caption{Splitting of $\Omega$ between finite volumes and finite elements.}
	\label{fig:triangles}
\end{figure} 
We now describe a rather general space discretization of the model problem. The core property we need is that the meshes of $\Omega_1$ and $\Omega_2$ share the same nodes on $\Gamma$ as shown in Figure \ref{fig:triangles}. Furthermore, we need that there is a specific set of unknowns associated with the interface nodes. Otherwise, we allow at this point for arbitraty meshes on both sides. 

Then, letting $\textbf{u}_I^{(1)}$ correspond to the unknowns on $\Omega_1$ and $\textbf{u}_{\Gamma}$ to the unknowns at the interface $\Gamma$, we can write a general discretization of the first equation in \eqref{model:eq} in a compact form as:

\begin{align}
	\label{eq1:fvm}
	\textbf{M}_1 \dot{\textbf{u}}_I^{(1)} + \textbf{M}_{I \Gamma}^{(1)} \dot{\textbf{u}}_{\Gamma} + \textbf{A}_1 \textbf{u}_I^{(1)} + \textbf{A}_{I \Gamma}^{(1)} \textbf{u}_{\Gamma} = \mathbf{0}.
\end{align}

On the other hand, a general discretization of the first equation in \eqref{model:eq} on $\Omega_2$ can be written as:

\begin{align}
\label{eq2:fvm}
\textbf{M}_2 \dot{\textbf{u}}_I^{(2)} + \textbf{M}_{I \Gamma}^{(2)} \dot{\textbf{u}}_{\Gamma} + \textbf{A}_2 \textbf{u}_I^{(2)} + \textbf{A}_{I \Gamma}^{(2)} \textbf{u}_{\Gamma} = \mathbf{0}.
\end{align}
where $\textbf{u}_I^{(2)}$ correspond to the unknowns on $\Omega_2$.

To close the system, we need an approximation of the normal derivatives at $\Gamma$. For the FVM on $\Omega_1$, we approximate the normal derivative with respect to $u_1$ using second order one-sided finite differences:

\begin{align}
\label{approxnormalderfvm}
-\lambda_1 \frac{\partial u_1}{\partial \textbf{n}_1} \approx \frac{\lambda_1}{2\Delta x} ( 4 u_{1,N} (t) - u_{1,N-1} (t) - 3 u_{\Gamma}).
\end{align}

On the other hand, let $\phi_j$ be a nodal FE basis function on $\Omega_2$ for a node on $\Gamma$ we observe that the normal derivative with respect to $u_2$ can be written as a linear functional using Green's formula \cite[pp. 3]{Toselli:2004qr}. Thus, the approximation of the normal derivative is given by

\begin{align}
\begin{split}
\label{approxnormalderfem}
& \lambda_2 \int_{\Gamma} \frac{\partial u_{2}}{\partial \textbf{n}_2} \phi_j dS = \lambda_2 \int_{\Omega_2} (\Delta u_{2} \phi_j + \nabla u_{2} \nabla \phi_j) d\textbf{x} \\
& = \alpha_2 \int_{\Omega_2}  \frac{d}{d t} u_{2} \phi_j + \lambda_2 \int_{\Omega_2} \nabla u_{2}  \nabla \phi_j d\textbf{x}.
\end{split}
\end{align}

Consequently, the equation

\begin{align}
\label{eq3:fvm}
\textbf{M}_{\Gamma \Gamma}^{(2)} \dot{\textbf{u}}_{\Gamma} + \textbf{M}_{\Gamma I}^{(2)} \dot{\textbf{u}}_{I}^{(2)} + \textbf{A}_{\Gamma \Gamma}^{(2)} \textbf{u}_{\Gamma} + \textbf{A}_{\Gamma I}^{(2)} \textbf{u}_{I}^{(2)} = - \textbf{M}_{\Gamma \Gamma}^{(1)} \dot{\textbf{u}}_{\Gamma} - \textbf{M}_{\Gamma I}^{(1)} \dot{\textbf{u}}_{I}^{(1)} - \textbf{A}_{\Gamma \Gamma}^{(1)} \textbf{u}_{\Gamma} - \textbf{A}_{\Gamma I}^{(1)} \textbf{u}_{I}^{(1)},
\end{align}
is a discrete version of the fourth equation in \eqref{model:eq} and completes the system \eqref{eq1:fvm}-\eqref{eq2:fvm}. Notice that the left hand side of \eqref{eq3:fvm} comes from \eqref{approxnormalderfem} and the right hand side from \eqref{approxnormalderfvm}. We can now write the coupled equations \eqref{eq1:fvm}, \eqref{eq2:fvm} and \eqref{eq3:fvm} as an ODE for the vector of unknowns $\textbf{u} = \left( \textbf{u}_I^{(1)}, \textbf{u}_I^{(2)}, \textbf{u}_{\Gamma} \right)^T$

\begin{align}
\tilde{\textbf{M}} \dot{\textbf{u}} + \tilde{\textbf{A}} \textbf{u} = \mathbf{0},
\end{align}
where

\begin{align}
\begin{split}
\tilde{\textbf{M}} = \left( \begin{array}{ccc}
\textbf{M}_1 & \mathbf{0} & \textbf{M}_{I \Gamma}^{(1)} \\
\mathbf{0} & \textbf{M}_2 & \textbf{M}_{I \Gamma}^{(2)} \\
\textbf{M}_{\Gamma I}^{(1)} & \textbf{M}_{\Gamma I}^{(2)} & \textbf{M}_{\Gamma \Gamma}^{(1)} + \textbf{M}_{\Gamma \Gamma}^{(2)} 
\end{array} \right), \  \  \tilde{\textbf{A}} = \left( \begin{array}{ccc}
\textbf{A}_1 & \mathbf{0} & \textbf{A}_{I \Gamma}^{(1)} \\
\mathbf{0} & \textbf{A}_2 & \textbf{A}_{I \Gamma}^{(2)} \\
\textbf{A}_{\Gamma I}^{(1)} & \textbf{A}_{\Gamma I}^{(2)} & \textbf{A}_{\Gamma \Gamma}^{(1)} + \textbf{A}_{\Gamma \Gamma}^{(2)} 
\end{array} \right). \nonumber
\end{split}
\end{align}

\subsection{Time Discretization}

Applying the implicit Euler method with time step $\Delta t$ to the system \eqref{eq3:fvm}, we get for the vector of unknowns $\textbf{u}^{n+1} = (\textbf{u}_I^{(1),n+1}, \textbf{u}_I^{(2),n+1}, \textbf{u}_{\Gamma}^{n+1})^T$

\begin{align}
\label{finallineq}
\textbf{A} \textbf{u}^{n+1} = \tilde{\textbf{M}} \textbf{u}^n,
\end{align}
where

\begin{align}
\textbf{A} = \tilde{\textbf{M}} + \Delta t \tilde{\textbf{A}} = \left( \begin{array}{ccc}
\textbf{M}_1 + \Delta t \textbf{A}_1 & \textbf{0} & \textbf{M}_{I \Gamma}^{(1)} + \Delta t \textbf{A}_{I \Gamma}^{(1)} \\
\textbf{0} & \textbf{M}_2 + \Delta t \textbf{A}_2 & \textbf{M}_{I \Gamma}^{(2)} + \Delta t \textbf{A}_{I \Gamma}^{(2)} \\
\textbf{M}_{\Gamma I}^{(1)} + \Delta t \textbf{A}_{\Gamma I}^{(1)} & \textbf{M}_{\Gamma I}^{(2)} + \Delta t \textbf{A}_{\Gamma I}^{(2)} & \textbf{M}_{\Gamma \Gamma} + \Delta t \textbf{A}_{\Gamma \Gamma}
\end{array} \right),   \nonumber
\end{align}
with $\textbf{M}_{\Gamma \Gamma} = \textbf{M}_{\Gamma \Gamma}^{(1)} + \textbf{M}_{\Gamma \Gamma}^{(2)}$ and $\textbf{A}_{\Gamma \Gamma} = \textbf{A}_{\Gamma \Gamma}^{(1)} + \textbf{A}_{\Gamma \Gamma}^{(2)}$.

\subsection{Dirichlet-Neumann Iteration}

We now employ a Dirichlet-Neumann iteration to solve the discrete system \eqref{finallineq}. This corresponds to alternately solving the discretized equations of the transmission problem \eqref{model:eq} on $\Omega_1$ with Dirichlet data on $\Gamma$ and the discretization of \eqref{model:eq} on $\Omega_2$ with Neumann data on $\Gamma$.

Therefore, from \eqref{finallineq} one gets for the $k$-th iteration the two equation systems 

\begin{align}
\label{fpi1}
(\textbf{M}_1 + \Delta t \textbf{A}_1) \textbf{u}_I^{(1),n+1,k+1} =  - (\textbf{M}_{I \Gamma}^{(1)} + \Delta t \textbf{A}_{I \Gamma}^{(1)}) \textbf{u}_{\Gamma}^{n+1,k} + \textbf{M}_1 \textbf{u}_I^{(1),n} + \textbf{M}_{I \Gamma}^{(1)} \textbf{u}_{\Gamma}^n,
\end{align}

\begin{align}
\label{fpi2}
\hat{\textbf{A}} \hat{\textbf{u}}^{k+1} = \hat{\textbf{M}} \textbf{u}^n - \textbf{b}^k,
\end{align}
to be solved in succession. Here,

\begin{align}
\hat{\textbf{A}} = \left( \begin{array}{cc}
\textbf{M}_2 + \Delta t \textbf{A}_2 & \textbf{M}_{I \Gamma}^{(2)} + \Delta t \textbf{A}_{I \Gamma}^{(2)} \\
\textbf{M}_{\Gamma I}^{(2)} + \Delta t \textbf{A}_{\Gamma I}^{(2)} & \textbf{M}_{\Gamma \Gamma}^{(2)} + \Delta t \textbf{A}_{\Gamma \Gamma}^{(2)}
\end{array} \right), \  \  \hat{\textbf{M}} = \left( \begin{array}{ccc}
\textbf{0} & \textbf{M}_2 & \textbf{M}_{I \Gamma}^{(2)} \\
\textbf{M}_{\Gamma I}^{(1)} & \textbf{M}_{\Gamma I}^{(2)} & \textbf{M}_{\Gamma \Gamma}
\end{array} \right), \nonumber
\end{align}
and

\begin{align}
\textbf{b}^k = \left( \begin{array}{c}
\textbf{0} \\
(\textbf{M}_{\Gamma I}^{(1)} + \Delta t \textbf{A}_{\Gamma I}^{(1)}) \textbf{u}_I^{(1),n+1,k+1} + (\textbf{M}_{\Gamma \Gamma}^{(1)} + \Delta t \textbf{A}_{\Gamma \Gamma}^{(1)}) \textbf{u}_{\Gamma}^{n+1,k}
\end{array} \right),
\end{align}

\begin{align}
\hat{\textbf{u}}^{k+1} = \left( \begin{array}{c}
\textbf{u}_I^{(2),n+1,k+1} \\
\textbf{u}_{\Gamma}^{n+1,k+1}
\end{array} \right), \nonumber
\end{align}
with some initial condition, here $\textbf{u}_{\Gamma}^{n+1,0} = \textbf{u}_{\Gamma}^{n}$. The iteration is terminated according to the standard criterion $\| \textbf{u}_{\Gamma}^{k+1} - \textbf{u}_{\Gamma}^k \| \leq \tau$ where $\tau$ is a user defined tolerance \cite{birken:15}. 

One way to analyze this method is to write it as a splitting method for \eqref{finallineq} and try to estimate the spectral radius of that iteration. However, the results obtained in this way are much too inaccurate. For that reason, we now rewrite \eqref{fpi1}-\eqref{fpi2} as an iteration for $\textbf{u}_{\Gamma}^{n+1}$ to restrict the size of the space to the dimension of $\textbf{u}_{\Gamma}$ which is much smaller. To this end, we isolate the term  $\textbf{u}_I^{(1),n+1,k+1}$ in \eqref{fpi1} and $\textbf{u}_I^{(2),n+1,k+1}$ in the first equation in \eqref{fpi2} and we insert the resulting expressions into the second equation in \eqref{fpi2}. Consequently, the iteration $\textbf{u}_{\Gamma}^{n+1,k+1} = \Sigma \textbf{u}_{\Gamma}^{n+1,k} + \psi^n$ is obtained with iteration matrix  

\begin{align}
\Sigma = - {\textbf{S}^{(2)}}^{-1} \textbf{S}^{(1)},
\end{align}
where

\begin{align}
\label{schur}
\textbf{S}^{(m)} = (\textbf{M}_{\Gamma \Gamma}^{(m)} + \Delta t \textbf{A}_{\Gamma \Gamma}^{(m)}) - (\textbf{M}_{\Gamma I}^{(m)} + \Delta t \textbf{A}_{\Gamma I}^{(m)}) (\textbf{M}_m + \Delta t \textbf{A}_m)^{-1} (\textbf{M}_{I \Gamma}^{(m)} + \Delta t \textbf{A}_{I \Gamma}^{(m)}), 
\end{align}
for $m=1,2$ and $\psi^n$ contains terms that depend only on the solutions at the previous time step. Notice that $\Sigma$ is a discrete version of the Steklov-Poincar\'e operator.  

Thus, the Dirichlet-Neumann iteration is a linear iteration and the rate of convergence is described by the spectral radius of the iteration matrix $\Sigma$.

\section{One-Dimensional Convergence Analysis}

The derivation so far was for a rather general discretization. In this section, we study the iteration matrix $\Sigma$ for a specific FVM-FEM discretization in 1D. We will give an exact formula for the convergence rates. The behaviour of the rates when approaching both the continuous case in time and space is also given.  

\begin{figure}[h]
	\centering
	\includegraphics[width=1\textwidth]{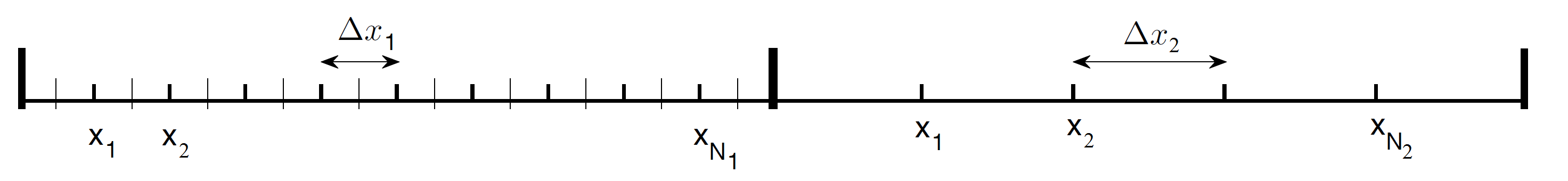}
	\caption{Grid cells over $\Omega_1$ and $\Omega_2$ for the finite volume discretization and the finite element discretization respectively.}
	\label{fig:fvmgrid}
\end{figure}

Specifically, we use $\Omega_1 = [-1,0]$, $\Omega_2 = [0,1]$. For the FVM discretization, we consider a primal grid, i.e, we discretize $\Omega_1$ into $N_1$ equal sized grid cells of size $\Delta x_1 = 1/(N_1+1)$, and define $x_i = i \Delta x_1$, so that $x_i$ is the center of the cell $i$, see figure \ref{fig:fvmgrid}. The edges of cell $i$ are then $x_{i-1/2}$ and $x_{i+1/2}$ and they form the corresponding dual grid. Moreover, we use the flux function

\begin{align}
F (u_L,u_R) = -\frac{\lambda_1}{\Delta x_1} (u_{1,i} - u_{1,i-1}), 
\end{align}
to approximate the flux, which results in a second order scheme. For the FEM discretization, we use the standard piecewise-linear polynomials as test functions. Here we discretize $\Omega_2$ into $N_2$ equal sized cells of size $\Delta x_2 = 1/(N_2+1)$. 

For the coupling between a compressible fluid and a structure, there would be a boundary layer in the fluid, meaning that the mesh would be very fine in direction normal to the boundary, implying $\Delta x_1 \ll \Delta x_2$. 

With $\textbf{e}_{m,j} = \left( \begin{array}{ccccccc}
0 & \cdots & 0 & 1 & 0 & \cdots & 0 \end{array} \right)^T \in \mathbb{R}^{N_m}$ where the only nonzero entry is located at the $j$-th position, the discretization matrices are given by

\begin{align}
\begin{split}
\textbf{A}_1 = \frac{\lambda_1}{\Delta x_1^2} \left( \begin{array}{cccc}
-2 & 1 & & 0 \\
1 & -2 & \ddots & \\
& \ddots & \ddots & 1 \\
0 & & 1 & -2 \\
\end{array} \right),  \  \  \textbf{A}_2 = \frac{\lambda_2}{\Delta x_2^2} \left( \begin{array}{cccc}
2 & -1 & & 0 \\
-1 & 2 & \ddots & \\
& \ddots & \ddots & -1 \\
0 & & -1 & 2 \\
\end{array} \right), \nonumber
\end{split}
\end{align} 

\begin{align}
\begin{split}
\textbf{M}_2 = \frac{\alpha_2}{6} \left( \begin{array}{cccc}
4 & 1 & & 0 \\
1 & 4 & \ddots & \\
& \ddots & \ddots & 1 \\
0 & & 1 & 4 \\
\end{array} \right), \  \  \textbf{A}_{\Gamma \Gamma}^{(1)} = \frac{3 \lambda_1}{2 \Delta x_1^2},  \  \  \textbf{A}_{\Gamma \Gamma}^{(2)} = \frac{\lambda_2}{\Delta x_2^2},  \  \  \textbf{M}_{\Gamma \Gamma}^{(2)} = \frac{2 \alpha_2}{6}, \nonumber
\end{split}
\end{align} 

\begin{align}
\begin{split}
\textbf{A}_{I \Gamma}^{(1)} = \frac{\lambda_1}{\Delta x_1^2} \textbf{e}_{1,N_1},  \  \  \textbf{A}_{I \Gamma}^{(2)} = -\frac{\lambda_2}{\Delta x_2^2} \textbf{e}_{2,1},  \  \  \textbf{M}_{I \Gamma}^{(2)} = \frac{\alpha_2}{6} \textbf{e}_{2,1}, \\
\textbf{A}_{\Gamma I}^{(1)} = \frac{\lambda_1}{2 \Delta x_1^2} (4 \textbf{e}_{1,N_1}^T - \textbf{e}_{1,N_1-1}^T),  \  \  \textbf{A}_{\Gamma I}^{(2)} = -\frac{\lambda_2}{\Delta x_2^2} \textbf{e}_{2,1}^T,  \  \  \textbf{M}_{\Gamma I}^{(2)} = \frac{\alpha_2}{6} \textbf{e}_{2,1}^T. \nonumber \\
\end{split}
\end{align}
where $\textbf{A}_m$, $\textbf{M}_m \in \mathbb{R}^{N_m \times N_m}$, $\textbf{A}_{I \Gamma}^{(m)}$, $\textbf{M}_{I \Gamma}^{(2)} \in \mathbb{R}^{N_m \times 1}$ and $\textbf{A}_{\Gamma I}^{(m)}$, $\textbf{M}_{\Gamma I}^{(2)} \in \mathbb{R}^{1 \times N_m}$ for $m=1,2$.

In this case, $\textbf{M}_1 = \alpha_1 \textbf{I}$, $\textbf{M}_{I \Gamma}^{(1)} = \textbf{M}_{\Gamma \Gamma}^{(1)} = \textbf{M}_{\Gamma I}^{(1)} =\mathbf{0}$. Thus, 

\begin{align}
\label{schur1_1D}
\textbf{S}^{(1)} = \Delta t \textbf{A}_{\Gamma \Gamma}^{(1)} - \Delta t^2 \textbf{A}_{\Gamma I}^{(1)} (\alpha_1 \textbf{I} - \Delta t \textbf{A}_1)^{-1} \textbf{A}_{I \Gamma}^{(1)},
\end{align}

\begin{align}
\label{schur2_1D}
\textbf{S}^{(2)} = (\textbf{M}_{\Gamma \Gamma}^{(2)} + \Delta t \textbf{A}_{\Gamma \Gamma}^{(2)}) - (\textbf{M}_{\Gamma I}^{(2)} + \Delta t \textbf{A}_{\Gamma I}^{(2)}) (\textbf{M}_2 + \Delta t \textbf{A}_2)^{-1} (\textbf{M}_{I \Gamma}^{(2)} + \Delta t \textbf{A}_{I \Gamma}^{(2)}).
\end{align}

Note that the iteration matrix $\Sigma$ is just a real number in this case and thus its spectral radius is its modulus. One computes $\textbf{S}^{(1)}$ and $\textbf{S}^{(2)}$ by inserting the corresponding matrices specified above in \eqref{schur1_1D} and \eqref{schur2_1D} obtaining

\begin{align}
\label{schur1D1}
\begin{split}
& \textbf{S}^{(1)} =  \Delta t \frac{3 \lambda_1}{2 \Delta x_1^2}  - \Delta t^2 \frac{\lambda_1^2}{2 \Delta x_1^4}  (4 \textbf{e}_{1,N_1}^T - \textbf{e}_{1,N_1-1}^T) (\alpha_1 \textbf{I} - \Delta t \textbf{A}_1)^{-1} \textbf{e}_{1,N_1} \\
& = \Delta t \frac{3\lambda_1}{2\Delta x_1^2}  - \Delta t^2 \frac{\lambda_1^2}{2\Delta x_1^4}  (4\alpha_{N_1N_1}^1 - \alpha_{N_1-1N_1}^1), \\
\end{split}
\end{align}

\begin{align}
\label{schur1D2}
\begin{split}
& \textbf{S}^{(2)} = \left( \frac{\alpha_2}{3} + \Delta t \frac{\lambda_2}{\Delta x_2^2} \right) - \left( \frac{\alpha_2}{6} - \Delta t \frac{\lambda_2}{\Delta x_2^2} \right)^2 \textbf{e}_{2,1}^T (\textbf{M}_2 + \Delta t \textbf{A}_2)^{-1} \textbf{e}_{2,1} \\
& = \left( \frac{\alpha_2}{3} + \Delta t \frac{\lambda_2}{\Delta x_2^2} \right) - \left( \frac{\alpha_2}{6} - \Delta t \frac{\lambda_2}{\Delta x_2^2} \right)^2 \alpha_{11}^2, \\
\end{split}
\end{align}
where $\alpha_{ij}^1$ represents the entries of the matrix $(\alpha_1 \textbf{I} - \Delta t \textbf{A}_1)^{-1}$ and $\alpha_{ij}^2$ the entries of $(\textbf{M}_2 + \Delta t \textbf{A}_2)^{-1}$ for $i,j=1,...,N_1$ and $i,j=1,...,N_2$ respectively. Observe that the matrices  $(\alpha_1 \textbf{I} - \Delta t \textbf{A}_1)^{-1}$ and $(\textbf{M}_2 + \Delta t \textbf{A}_2)$ are tridiagonal Toeplitz matrices but their inverses are full matrices. The computation of the exact inverses is based on a recursive formula which runs over the entries \cite{Fonseca:2001} and consequently, it is non trivial to compute $\alpha_{N_1N_1}^1$, $\alpha_{N_1-1N_1}^1$ and $\alpha_{11}^2$ this way. 

Due to these difficulties, we propose to rewrite the matrices $(\alpha_1 \textbf{I} - \Delta t \textbf{A}_1)^{-1}$ and $(\textbf{M}_2 + \Delta t \textbf{A}_2)^{-1}$ in terms of their eigendecomposition:

\begin{align}
\begin{split}
(\alpha_1 \textbf{I} - \Delta t \textbf{A}_1)^{-1} = \left[ \mbox{tridiag} \left( -\frac{\lambda_1 \Delta t}{ \Delta x_1^2}, \frac{ \alpha_1 \Delta x_1^2 + 2 \lambda_1 \Delta t}{ \Delta x_1^2}, -\frac{\lambda_1 \Delta t}{ \Delta x_1^2} \right) \right]^{-1} = \textbf{V}_{N_1} \Lambda_1^{-1} \textbf{V}_{N_1},
\end{split}
\end{align}

\begin{align}
\begin{split}
(\textbf{M}_2 + \Delta t \textbf{A}_2)^{-1} = \left[ \mbox{tridiag} \left( \frac{\alpha_2 \Delta x_2^2 - 6 \lambda_2 \Delta t}{6 \Delta x_2^2}, \frac{2 \alpha_2 \Delta x_2^2 + 6 \lambda_2 \Delta t}{3 \Delta x_2^2}, \frac{\alpha_2 \Delta x_2^2 - 6 \lambda_2 \Delta t}{6 \Delta x_2^2} \right) \right]^{-1} = \textbf{V}_{N_2} \Lambda_2^{-1} \textbf{V}_{N_2},
\end{split}
\end{align}
where the matrix $\textbf{V}_N$ has the eigenvectors of any symmetric tridiagonal Toeplitz matrix of dimension $N$ as columns. The entries of $\textbf{V}_{N_1}$ and $\textbf{V}_{N_2}$ are not dependent on the entries of $\alpha_1 \textbf{I} - \Delta t \textbf{A}_1$ or $\textbf{M}_2 + \Delta t \textbf{A}_2$ due to their symmetry. Moreover, the matrices $\Lambda_1$ and $\Lambda_2$ are diagonal matrices having the eigenvalues of $\alpha_1 \textbf{I} - \Delta t \textbf{A}_1$ or $\textbf{M}_2 + \Delta t \textbf{A}_2$ as entries respectively. These are known and given e.g. in \cite[pp. 514-516]{Meyer:2000}:

\begin{align}
\begin{split}
& v_{ij}^1 = \frac{1}{\sqrt{\sum_{k=1}^{N_1} \sin^2 \left( \frac{k \pi}{N_1+1} \right)}} \sin \left( \frac{ij \pi}{N_1+1} \right) \  \  \mbox{for} \  \  i,j=1,...,N_1, \\
& v_{ij}^2 = \frac{1}{\sqrt{\sum_{k=1}^{N_2} \sin^2 \left( \frac{k \pi}{N_2+1} \right)}} \sin \left( \frac{ij \pi}{N_2+1} \right) \  \  \mbox{for} \  \  i,j=1,...,N_2, \\
& \mu_{1,j} = \frac{1}{\Delta x_1^2} \left( \alpha_1 \Delta x_1^2 + 2 \lambda_1 \Delta t - 2 \lambda_1 \Delta t \cos \left( \frac{j \pi}{N_1+1} \right) \right) \  \  \mbox{for} \  \  j = 1,...,N_1, \\
& \mu_{2,j} = \frac{1}{3 \Delta x_2^2} \left( 2 \alpha_2 \Delta x_2^2 + 6 \lambda_2 \Delta t + (\alpha_2 \Delta x_2^2 - 6 \lambda_2 \Delta t) \cos \left( \frac{j \pi}{N_2+1} \right) \right) \  \  \mbox{for} \  \  j = 1,...,N_2.
\end{split}
\end{align}

The entries $\alpha_{N_1N_1}^1$, $\alpha_{N_1-1N_1}^1$ and $\alpha_{11}^2$ of the matrices $(\alpha_1 \textbf{I} - \Delta t \textbf{A}_1)^{-1}$ and $(\textbf{M}_2 + \Delta t \textbf{A}_2)^{-1}$, respectively, are now computed through their eigendecomposition resulting in 

\begin{align}
\begin{split}
\label{alphaN1N}
\alpha_{N_1-1 N_1}^1 = \frac{\Delta x_1^2 s_0}{\sum_{i=1}^{N_1} \sin^2 (i \pi \Delta x_1)},
\end{split}
\end{align}

\begin{align}
\label{alphaNN}
\begin{split}
\alpha_{N_1N_1}^1 = \frac{\Delta x_1^2 s_1}{\sum_{i=1}^{N_1} \sin^2 (i \pi \Delta x_1)},
\end{split}
\end{align}

\begin{align}
\label{alpha11}
\begin{split}
\alpha_{11}^2 = \frac{3 \Delta x_2^2 s_2}{\sum_{i=1}^{N_2} \sin^2 (i \pi \Delta x_2)},
\end{split}
\end{align}
with

\begin{align}
s_0 = \sum_{i=1}^{N_1} \frac{ \sin (i \pi \Delta x_1) \sin (2 i \pi \Delta x_1)}{\alpha_1 \Delta x_1^2 + 2 \lambda_1 \Delta t (1 - \cos (i \pi \Delta x_1))},
\end{align}

\begin{align}
s_1 = \sum_{i=1}^{N_1} \frac{ \sin^2 (i \pi \Delta x_1)}{\alpha_1 \Delta x_1^2 + 2 \lambda_1 \Delta t (1 - \cos (i \pi \Delta x_1))},
\end{align}

\begin{align}
\label{sumden1D}
s_2 = \sum_{i=1}^{N_2} \frac{ \sin^2 (i \pi \Delta x_2)}{2 \alpha_2 \Delta x_2^2 + 6 \lambda_2 \Delta t + (\alpha_2 \Delta x_2^2 - 6 \lambda_2 \Delta t) \cos (i \pi \Delta x_2)}.
\end{align}

Now, inserting \eqref{alphaN1N}, \eqref{alphaNN} and \eqref{alpha11} into \eqref{schur1D1} and \eqref{schur1D2} we get for $\textbf{S}^{(1)}$ and $\textbf{S}^{(2)}$:

\begin{align}
\begin{split}
\textbf{S}^{(1)} = \frac{3 \lambda_1 \Delta t}{2 \Delta x_1^2} -  \frac{\lambda_1^2 \Delta t^2}{2 \Delta x_1^2} \frac{4 s_1 - s_0}{\sum_{i=1}^{N_1} \sin^2 (i \pi \Delta x_1)},  \\
\end{split}
\end{align}

\begin{align}
\begin{split}
\label{S2:finaleqfem}
\textbf{S}^{(2)} = \left( \frac{\alpha_2 \Delta x_2^2 + 3 \lambda_2 \Delta t}{3 \Delta x_2^2} \right) -  \frac{(\alpha_2 \Delta x_2^2 - 6 \lambda_2 \Delta t)^2}{12 \Delta x_2^2} \frac{s_2}{\sum_{i=1}^{N_2} \sin^2 (i \pi \Delta x_2)}.  \\
\end{split}
\end{align}

With this we obtain an explicit formula for the spectral radius of the iteration matrix $\Sigma$ as a function of $\Delta x_1$, $\Delta x_2$ and $\Delta t$:

\begin{align}
\label{specrad1D}
\begin{split}
& \rho (\Sigma) = | \Sigma | = | {\textbf{S}^{(2)}}^{-1} \textbf{S}^{(1)} | \\
& = \left(  \frac{\alpha_2 \Delta x_2^2 + 3 \lambda_2 \Delta t}{3 \Delta x_2^2}  -  \frac{(\alpha_2 \Delta x_2^2 - 6 \lambda_2 \Delta t)^2}{12 \Delta x_2^2} \frac{s_2}{\sum_{i=1}^{N_2} \sin^2 (i \pi \Delta x_2)} \right)^{-1} \\
& \cdot \left( \frac{ 3\lambda_1 \Delta t}{ 2\Delta x_1^2}  -  \frac{\lambda_1^2 \Delta t^2}{2\Delta x_1^2} \frac{4s_1 - s_0}{\sum_{i=1}^{N_1} \sin^2 (i \pi \Delta x_1)} \right). 
\end{split}
\end{align}

To simplify this, the finite sums $\sum_{i=1}^{N_1} \sin^2 (i \pi \Delta x_1)$ and $\sum_{i=1}^{N_2} \sin^2 (i \pi \Delta x_2)$ can be computed. We first rewrite the sum of squared sinus terms into a sum of cosinus terms using the identity $\sin^2 (x/2) = (1-\cos(x))/2$. Then, the resulting sum can be converted into a geometric sum using Euler's formula. We thus obtain after some calculations:

\begin{align}
\label{sum1D}
\begin{split}
\sum_{j=1}^{N_1} \sin^2 (j \pi \Delta x_1) = \frac{1-\Delta x_1}{2\Delta x_1} - \frac{1}{2} \sum_{j=1}^{N_1} \cos (2 j \pi \Delta x_1) = \frac{1}{2 \Delta x_1},
\end{split}
\end{align} 

\begin{align}
\label{sum1D2}
\sum_{j=1}^{N_2} \sin^2 (j \pi \Delta x_2) = \frac{1}{2 \Delta x_2}.
\end{align}

Inserting \eqref{sum1D} and \eqref{sum1D2} into \eqref{specrad1D} we get after some manipulations

\begin{align}
\label{formula1D}
\begin{split}
| \Sigma | = \frac{3 \Delta x_2^2( 3\lambda_1 \Delta t - 2 \lambda_1^2 \Delta x_1 \Delta t^2 (4 s_1 - s_0))}{\Delta x_1^2 (2 (\alpha_2 \Delta x_2^2 + 3 \lambda_2 \Delta t)  - \Delta x_2 (\alpha_2 \Delta x_2^2 - 6 \lambda_2 \Delta t)^2 s_2)}.
\end{split}
\end{align}

We could not find a way of simplifying the finite sum \eqref{sumden1D} because $\Delta x_2$ depends on $N_2$ (i.e., $\Delta x_2 = 1/(N_2+1)$). However, \eqref{formula1D} is a computable formula that gives exactly the convergence rates of the Dirichlet-Neumann iteration for given $\Delta t$, $\Delta x_m$, $\alpha_m$ and $\lambda_m$, $m=1,2$.  

We are now interested in the asymptotics of \eqref{formula1D} for $\Delta t \rightarrow 0$ and $\Delta x_1 \rightarrow 0$ with $\Delta x_2 = r \cdot \Delta x_1$ where $r:= \Delta x_2 / \Delta x_1$ is a fixed aspect ratio. This is motivated by the assumption that the resolution in the fluid in direction tangential to the wall would be similar to the resolution in the structure. We obtain: 

\begin{align}
\label{limitdeltat1D}
\lim_{\Delta t \rightarrow 0} |\Sigma| = \frac{3 \Delta x_2^2 \cdot 0}{\Delta x_1^2 \left(2 \alpha_2 \Delta x_2^2 - \alpha_2 \Delta x_2^3 \sum_{i=1}^{N_2} \frac{3 \sin^2 (i \pi \Delta x_2)}{2 + \cos(i \pi \Delta x_2)} \right)} = 0.
\end{align}

\begin{align}
\label{limitdeltax1D}
\begin{split}
& \lim_{\Delta x_1 \rightarrow 0} |\Sigma| = \lim_{\Delta x_1 \rightarrow 0} \frac{3 r^2 ( 3\lambda_1 \Delta t - 2 \lambda_1^2 \Delta x_1 \Delta t^2 (4 s_1 - s_0))}{ 2 (\alpha_2 r^2 \Delta x_1^2 + 3 \lambda_2 \Delta t)  - r \Delta x_1 (\alpha_2 r^2 \Delta x_1^2 - 6 \lambda_2 \Delta t)^2 s_2} \\
& = \lim_{\Delta x_1 \rightarrow 0} \frac{9 \lambda_1 r^2 \Delta t - 6 \lambda_1 r^2 \Delta x_1 \Delta t \left( \sum_{i=1}^{N_1} \frac{\sin^2 (i \pi \Delta x_1) (2 - \cos(i\pi \Delta x_1))}{1 - \cos (i \pi \Delta x_1)} \right) }{6 \lambda_2 \Delta t - 6 \lambda_2 r \Delta t \Delta x_1 \left( \sum_{i=1}^{N_2} \frac{\sin^2 (i \pi r \Delta x_1)}{1 - \cos(i \pi r \Delta x_1)} \right)} \\
& = \frac{\lambda_1}{\lambda_2} \lim_{\Delta x_1 \rightarrow 0} \frac{3 r^2 - 2 r^2 \Delta x_1 \sum_{i=1}^{N_1} (1 + \cos (i \pi \Delta x_1)) (2 - \cos(i \pi \Delta x_1))}{2 - 2 r \Delta x_1 \sum_{i=1}^{N_2} (1 + \cos(i \pi r \Delta x_1))} \\
& = \frac{\lambda_1}{\lambda_2} \lim_{\Delta x_1 \rightarrow 0} \frac{3 r^2 - 2 r^2 \Delta x_1 \left( \sum_{i=1}^{N_1} 2 + \sum_{i=1}^{N_1} \cos(i \pi \Delta x_1 ) - \sum_{i=1}^{N_1} \cos^2 (i \pi \Delta x_1)\right)}{2 - 2r \Delta x_1 \left( \sum_{i=1}^{N_2} 1 + \sum_{i=1}^{N_2} cos (i \pi r \Delta x_1) \right)}.
\end{split}
\end{align}

To simplify \eqref{limitdeltax1D}, it is well known that the finite sums $\sum_{i=1}^{N_1} \cos(i \pi \Delta x_1)$, $\sum_{i=1}^{N_2} \cos(i \pi r \Delta x_1)$ and $\sum_{i=1}^{N_1} \cos^2 (i \pi \Delta x_1)$ can be computed by using Euler's formula to convert them into geometric sums. We thus obtain after some calculations:

\begin{align}
\label{newsum2}
\sum_{j=1}^{N_2} \cos(j \pi r \Delta x_1) = Re \left( \sum_{j=1}^{N_2} e^{ij \pi r \Delta x_1} \right) = Re \left( \frac{e^{i \pi r \Delta x_1} (1 - e^{i N_2 \pi r \Delta x_1})}{1 - e^{i \pi r \Delta x_1}} \right) = 0.
\end{align}

In order to compute the third sum, we rewrite the sum of squared cosinus terms into a sum of sinus terms using the identity $\cos^2 (x/2) = (1 + \cos(x))/2$ and the apply the same technique:

\begin{align}
\begin{split}
\label{newsum3}
\sum_{j=1}^{N_1} \cos^2 (j \pi \Delta x_1) = \frac{1-\Delta x_1}{2\Delta x_1} + \frac{1}{2} \sum_{j=1}^{N_1} \cos (2 j \pi \Delta x_1)  = \frac{1 - 2\Delta x_1}{2 \Delta x_1}.
\end{split}
\end{align}  

Inserting \eqref{newsum2} and \eqref{newsum3} into \eqref{limitdeltax1D} we get 

\begin{align}
\begin{split}
\label{finallimitdeltax1D}
& \lim_{\Delta x_1 \rightarrow 0} | \Sigma | = \frac{\lambda_1}{\lambda_2} \lim_{\Delta x_1 \rightarrow 0} \frac{3 r^2 - 2 r^2 \Delta x_1 \left( \frac{2(1 - \Delta x_1)}{\Delta x_1} - \frac{1 - 2\Delta x_1}{2\Delta x_1} \right)}{2 - 2r \Delta x_1 \left( \frac{1 - r \Delta x_1}{r \Delta x_1} \right)} \\
& = \frac{\lambda_1}{\lambda_2} \lim_{\Delta x_1 \rightarrow 0} \frac{2 r^2 \Delta x_1}{2 r \Delta x_1}= \frac{\lambda_1}{\lambda_2} r =: \delta_r.
\end{split}
\end{align}

From the result obtained in \eqref{limitdeltat1D} we can conclude that the convergence rate goes to zero when the time step decreases and therefore, the iteration will be fast for $\Delta t$ small and can always be made to converge by decreasing $\Delta t$. This is consistent with the behavior of the cooling of a flat plat and the flanged shaft presented earlier in figures \ref{fig:coolingsystemsrates}a and \ref{fig:coolingsystemsrates}b. 

On the other hand, from the spatial asymptotics \eqref{finallimitdeltax1D} we can observe that strong jumps in the thermal conductivities of the materials placed in $\Omega_1$ and $\Omega_2$ will imply fast convergence. This is often the case when modelling thermal fluid structure interaction, since fluids typically have lower thermal conductivities than structures.

Finally, the aspect ratio $r$ also influences the behavior of the fixed point iteration, i.e, the rates will become smaller the higher the aspect ratio, e.g. the higher the Reynolds number in the fluid. This phenomenon is not unknown for PDE discretizations and is referred to as geometric stiffness. As is the case here, refining the mesh to reduce the aspect ratio would lead to faster convergence of the iterative method. 

Before presenting numerical results we want to show the results obtained for different space discretization combinations with the same constant mesh width on both subdomains.

\section{Extension of the Analysis}

In this section we want to extend the results presented in the previous section by reviewing similar analysis for other choices of space discretizations. In particular, FEM-FEM coupling and 2D FVM-FEM with $r=1$.

Firstly, when one uses a linear FEM discretization in 1D and the same mesh width on both subdomains (i.e, $r=1$) and applies the same analysis as in the previous section, the corresponding limits for the spectral radius of the iteration matrix $\Sigma$ are given by \cite{Mongelic:16,Monge:16}:

\begin{align}
\label{limitFEMFEMdt}
\lim_{\Delta t \rightarrow 0} \rho (\Sigma) = \frac{\alpha_1}{\alpha_2},
\end{align}

\begin{align}
\label{limitFEMFEMdx}
\lim_{\Delta x \rightarrow 0} \rho (\Sigma) = \frac{\lambda_1}{\lambda_2}.
\end{align}

When we compare these with the asymptotics obtained with FVM-FEM discretizations \eqref{limitdeltat1D}-\eqref{finallimitdeltax1D}, we observe that while the spatial limit is the same, the temporal limit does not match. This arises from differences in the matrix $\textbf{S}^{(1)}$ in \eqref{schur}. In the FEM-FEM context, the matrices $\textbf{S}^{(1)}$ and $\textbf{S}^{(2)}$ lead to the same expression with only different material coefficients ($\alpha_1$, $\alpha_2$, $\lambda_1$, $\lambda_2$). Because of this, the limits of $\rho(\Sigma)$ are quotients of those coefficients. However, the situation is different in the FVM-FEM context. There, the matrix $\textbf{S}^{(1)}$ in \eqref{schur1_1D} is missing several mass matrices if we compare it with $\textbf{S}^{(2)}$ in \eqref{schur2_1D}. This unsymmetry between $\textbf{S}^{(1)}$ and $\textbf{S}^{(2)}$ causes that the limit of $\rho(\Sigma)$ when $\Delta t \rightarrow 0$ is not balanced between the numerator and the denominator, resulting in 0. 

This implies that opposed to the FVM-FEM case, where convergence can always be achieved by decreasing the time step, that for an FEM-FEM coupling, a situation can occur where $\alpha_1/\alpha_2>\lambda_1/\lambda_2$ and therefore, a decrease in time step can cause divergence. This is for example the case for an air-water coupling \cite{Mongelic:16}. 

Secondly, for an aspect ratio of $r=1$, we were able to extend the 1D results for both FVM-FEM and FEM-FEM to 2D in the following sense (see \cite{Birken:2016,Mongelic:16}). In 2D, the iteration matrix $\Sigma$ is not easy to compute for several reasons. First of all, the matrices $\textbf{M}_1 + \Delta t \textbf{A}_1$ and $\textbf{M}_2 + \Delta t \textbf{A}_2$ are sparse block tridiagonal matrices, and consequently their inverses are not straight forward to compute. Moreover, the diagonal blocks of the same matrices are tridiagonal but their inverses are full matrices.  

Due to these difficulties, we approximated the strictly diagonally dominant matrices $\textbf{M}_1 + \Delta t \textbf{A}_1$ and $\textbf{M}_2 + \Delta t \textbf{A}_2$ by their diagonal. Thus, we obtained an estimate of the spectral radius of the iteration matrix $\Sigma$.
% \begin{align}
% \rho (\Sigma) \approx \frac{(5 \alpha_2 \Delta x^2 + 24 \lambda_2 \Delta t) (2 (5 \alpha_1 \Delta x^2 + 24 \lambda_1 \Delta t)^2 - (\alpha_1 \Delta x^2 + 12 \lambda_1 \Delta t)^2)}{(5 \alpha_1 \Delta x^2 + 24 \lambda_1 \Delta t) (2 (5 \alpha_2 \Delta x^2 + 24 \lambda_2 \Delta t)^2 - (\alpha_2 \Delta x^2 + 12 \lambda_2 \Delta t)^2)}.
% \end{align}
This estimator tends to the exact same limits as for the 1D case for both combination of discretizations. 

We did not find a way to further extend these results to the high aspect ratio case. However, we will show now by numerical experiments that already the 1D formula \eqref{formula1D} is a good estimator for convergence rates in 2D.

\section{Numerical Results}

We now present numerical experiments designed to illustrate the validity of the theoretical results of the previous sections. Firstly, we will confirm that the theoretical formula $|\Sigma|$ in \eqref{formula1D} predicts the convergence rates in the 1D case. Secondly, we will show the validity of \eqref{formula1D} as an estimator for the rates in the 2D case, we will also show that the theoretical asymptotics deduced in \eqref{limitdeltat1D} and \eqref{finallimitdeltax1D} match with the numerical experiments. Finally, we illustrate the validity of \eqref{formula1D} as an estimator for the non linear thermal FSI test cases introduced in section 3.

\subsection{1D FVM-FEM Results}

We first compare the semidiscrete estimator $\beta$ in \eqref{semidiscretelimit} with the discrete formula $|\Sigma|$ in 1D in \eqref{formula1D} and experimental convergences rates. The latter are computed with respect to a reference solution $u_{ref}$ over the whole domain $\Omega$.    

\begin{figure}[h!]
	\centering
	\subfigure[$\Delta t / \Delta x^2 \ll 1$]{\includegraphics[width=6cm]{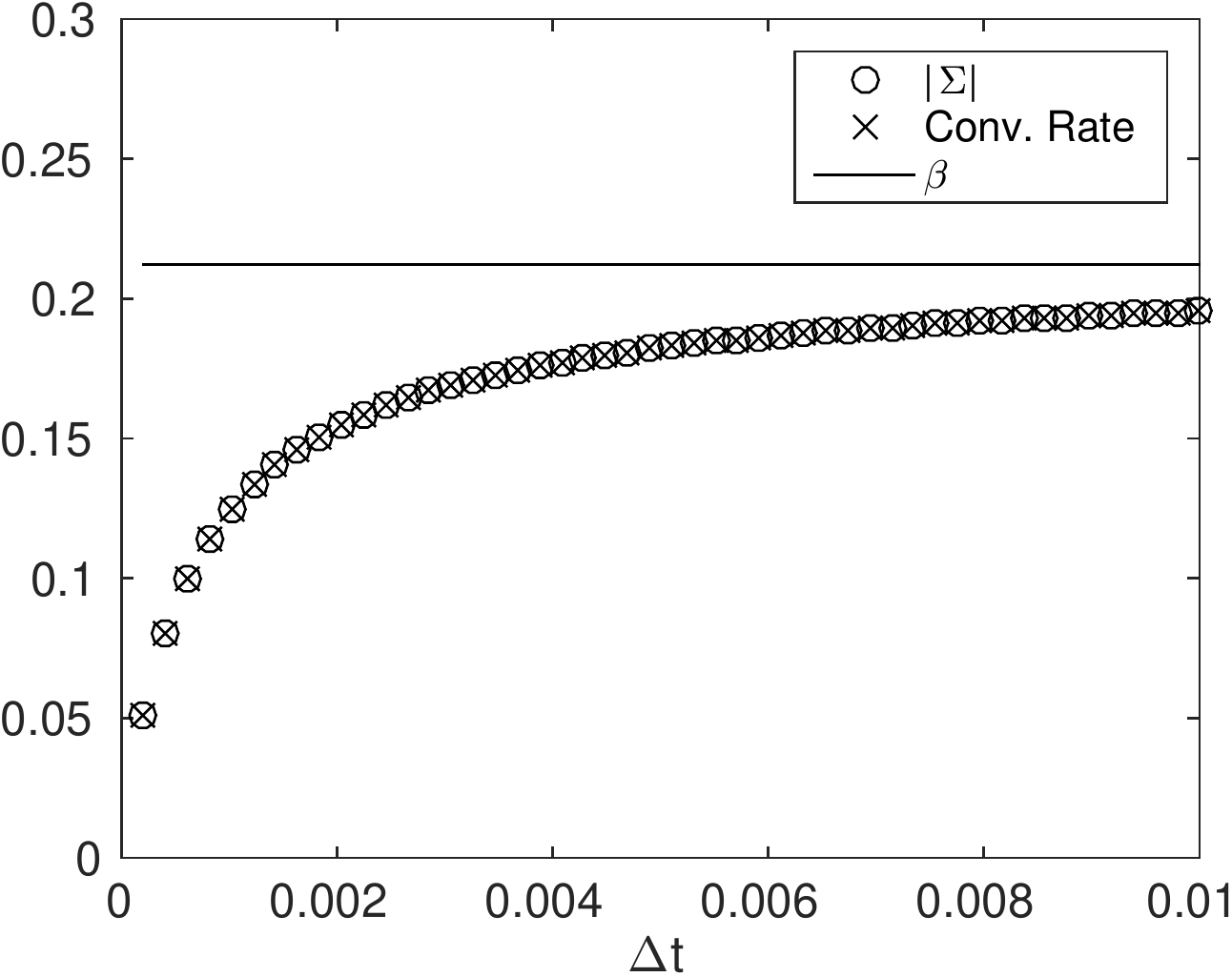}} \hfill
	\subfigure[$\Delta t / \Delta x^2 \gg 1$]{\includegraphics[width=6cm]{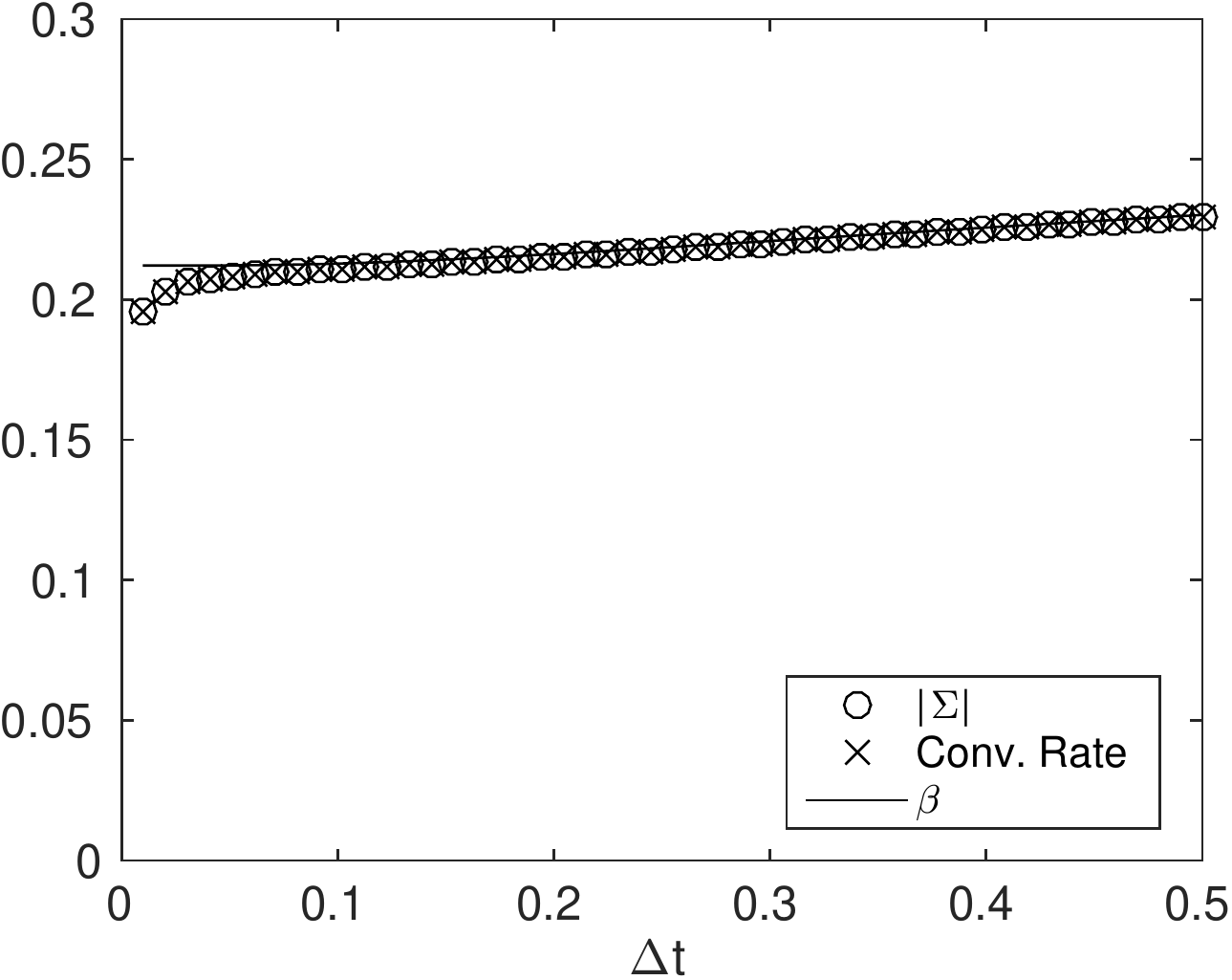}}
	\caption{Semidiscrete estimator $\beta$, exact rate $\Sigma$ and numerical rates over $\Delta t$ in 1D. $D_1 = 1$, $D_2 = 0.5$, $\lambda_1 = 0.3$ and $\lambda_2 = 1$, $\Delta x=1/20$. On the left: $\Delta t = 1e-2/50, 2 \cdot 1e-2/50,..., 50 \cdot 1e-2/50$. On the right: $\Delta t = 1e-2, 2 \cdot 1e-2,..., 50 \cdot 1e-2$.}
	\label{fig:disvssemi1D}
\end{figure}

Figure \ref{fig:disvssemi1D} shows a comparison between $\beta$ and $|\Sigma|$ for $r=1$, $\Delta x=1/20$ and variable $\Delta t$. On the left we plot $\beta$, $|\Sigma|$ and the experimental convergence rates with $\Delta t / \Delta x^2 \ll 1$ and on the right we plot the same but with $\Delta t / \Delta x^2 \gg 1$. As can be seen, the experimental convergence rate matches exactly with the exact formula \eqref{formula1D}. Observe that $\beta$ is almost constant and represents the lower branch in figure \ref{fig:semidiscretevsdeltat}. To arrive at the jump we would have to choose $\Delta t \gg 1$. We can conclude that the formulas for the convergence rates in 1D presented in the previous section are minimally better than the semidiscrete one proposed in \cite{Henshaw:2009} when $\Delta t / \Delta x^2 \gg 1$. In the, less relevant case, $\Delta t / \Delta x^2 \ll 1$ our formula also predicts the rates accurately, while the semidiscrete estimator deviates.

We now want to illustrate how $|\Sigma|$ in \eqref{formula1D} gives the convergence rates and tends to the limits computed previously in \eqref{limitdeltat1D} and \eqref{finallimitdeltax1D}. To this end, we present two real data examples. We consider here the thermal interaction between air at $273K$ with steel at $900K$ and water at $283K$ with steel at $900K$. Physical properties of the materials and resulting asymptotics for these two cases are shown in table~\ref{tab:physicalcoeff} and \ref{tab:asymptotics} respectively. 

\begin{table}
\caption{Physical properties of the materials. $\lambda$ is the thermal conductivity, $\rho$ the density, $c_p$ the specific heat capacity and $\alpha = \rho c_p$.}
\label{tab:physicalcoeff}     

\begin{tabular}{lllll}
\hline\noalign{\smallskip}
\textbf{Material} & $\lambda$ (W/mK) & $\rho$ (kg/$\mbox{m}^3$) & $c_p$ (J/kgK) & $\alpha$ (J/K $\mbox{m}^3$) \\
\noalign{\smallskip}\hline\noalign{\smallskip}
\textbf{Air} & 0.0243 & 1.293 & 1005 & 1299.5 \\
\textbf{Water} & 0.58 & 999.7 & 4192.1 & 4.1908e6  \\
\textbf{Steel} & 48.9 & 7836 & 443 & 3471348  \\
\noalign{\smallskip}\hline
\end{tabular}  
\end{table}

\begin{table}
\caption{Temporal and spatial asymptotics of \eqref{formula1D} for the thermal interaction of air at $273K$ with steel at $900K$, water at $283K$ with steel and air with water.}
\label{tab:asymptotics} 

\begin{tabular}{lll}
\hline\noalign{\smallskip}
\textbf{Case} & $\Delta t \rightarrow 0$ & $\Delta x \rightarrow 0$ \\
\noalign{\smallskip}\hline\noalign{\smallskip}
\textbf{Air-Steel} & 0 & 4.9693e-4 \cdot r \\
\textbf{Water-Steel} & 0 & 0.0119 \cdot r \\
\textbf{Air-Water} & 0 & 0.0419 \cdot r \\
\noalign{\smallskip}\hline
\end{tabular}  
\end{table}

\begin{figure}[h!]
	\centering
	\subfigure[$\Delta t = 40/39, 2 \cdot 40/39,..., 39 \cdot 40/39$, $\Delta x_1 = 1/1100$ and $r = 100$ (top curves) or $r=1$ (bottom curves).]{\includegraphics[width=5.9cm]{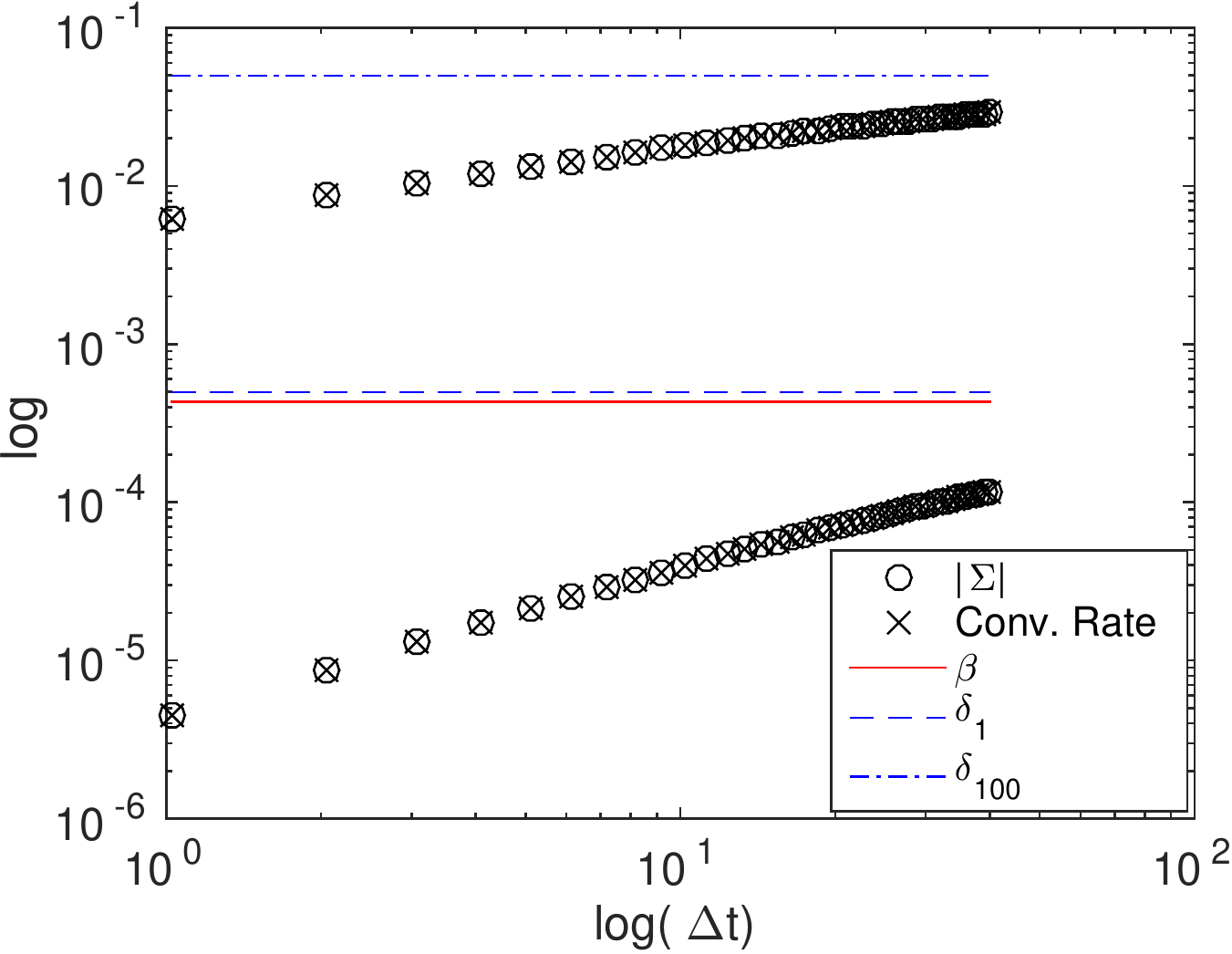}} \hfill
	\hspace{0.1cm}
	\subfigure[$\Delta x_1 = 1/3, 1/4,..., 1/50$, $\Delta t = 10$ and $r=100$ (top curves) or $r=1$ (bottom curves).]{\includegraphics[width=5.9cm]{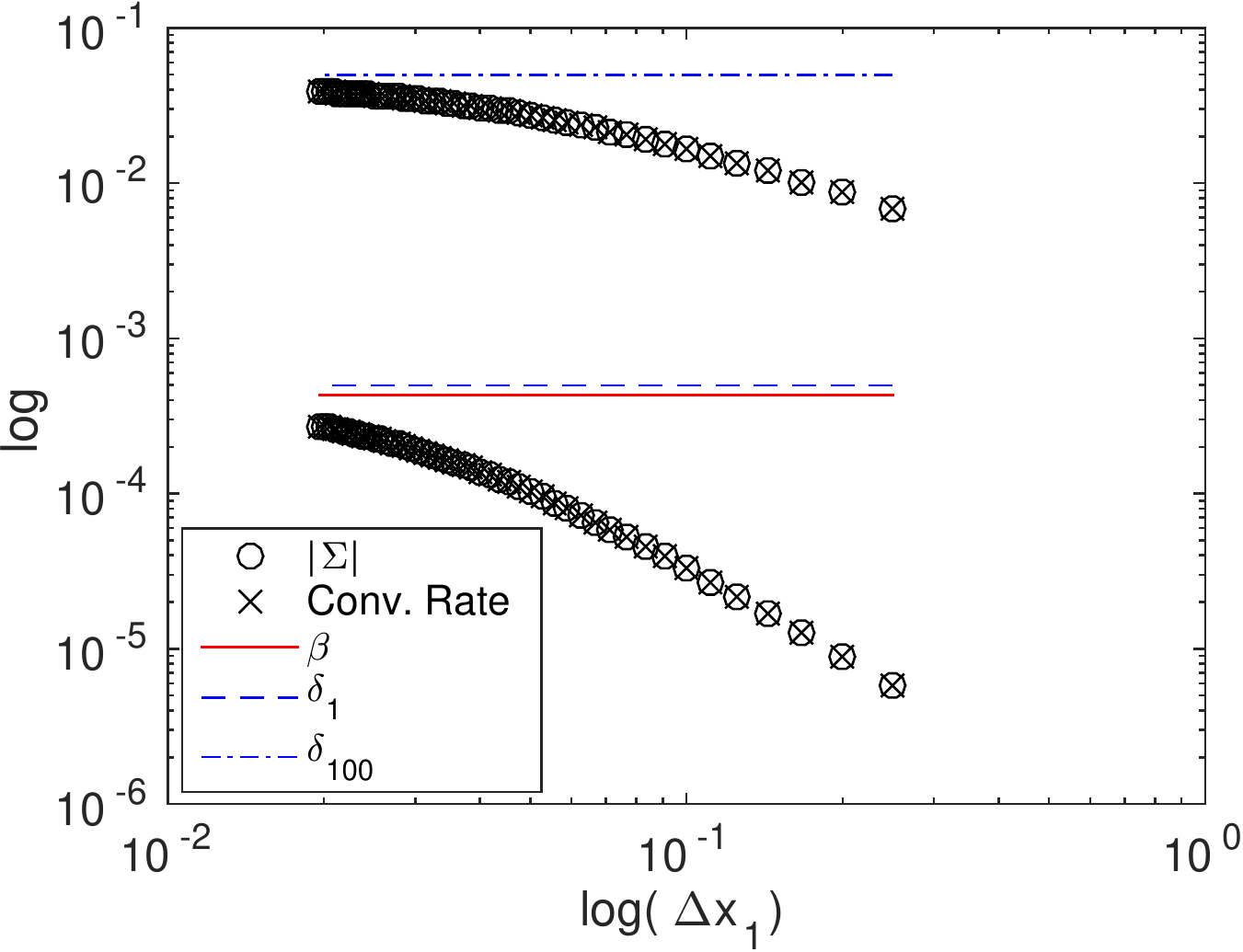}}
	\caption{Air-Steel thermal interaction with respect $\Delta t$ on the left and $\Delta x_1$ on the right in 1D.}
	\label{fig:air_steel}
\end{figure}

\begin{figure}[h!]
	\centering
	\subfigure[$\Delta t = 1/39, 2 \cdot 1/39,..., 39 \cdot 1/39$, $\Delta x_1 = 1/1100$ and $r = 1$ (top curves) or $r=0.01$ (bottom curves).]{\includegraphics[width=5.9cm]{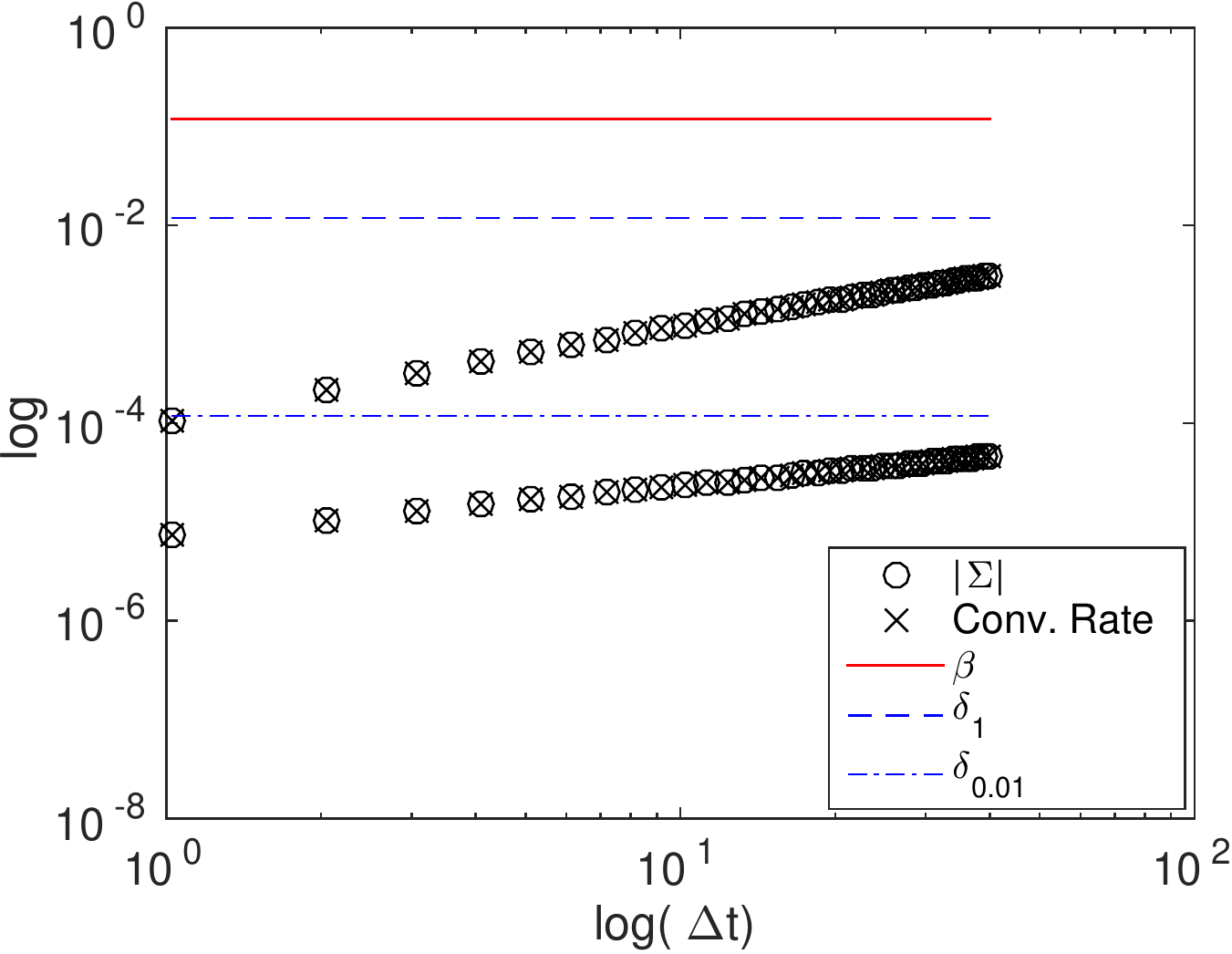}} \hfill
	\hspace{0.1cm}
	\subfigure[$\Delta x_1 = 1/3, 1/4,..., 1/50$, $\Delta t = 10$ and $r=1$ (top curves) or $r=0.01$ (bottom curves).]{\includegraphics[width=5.9cm]{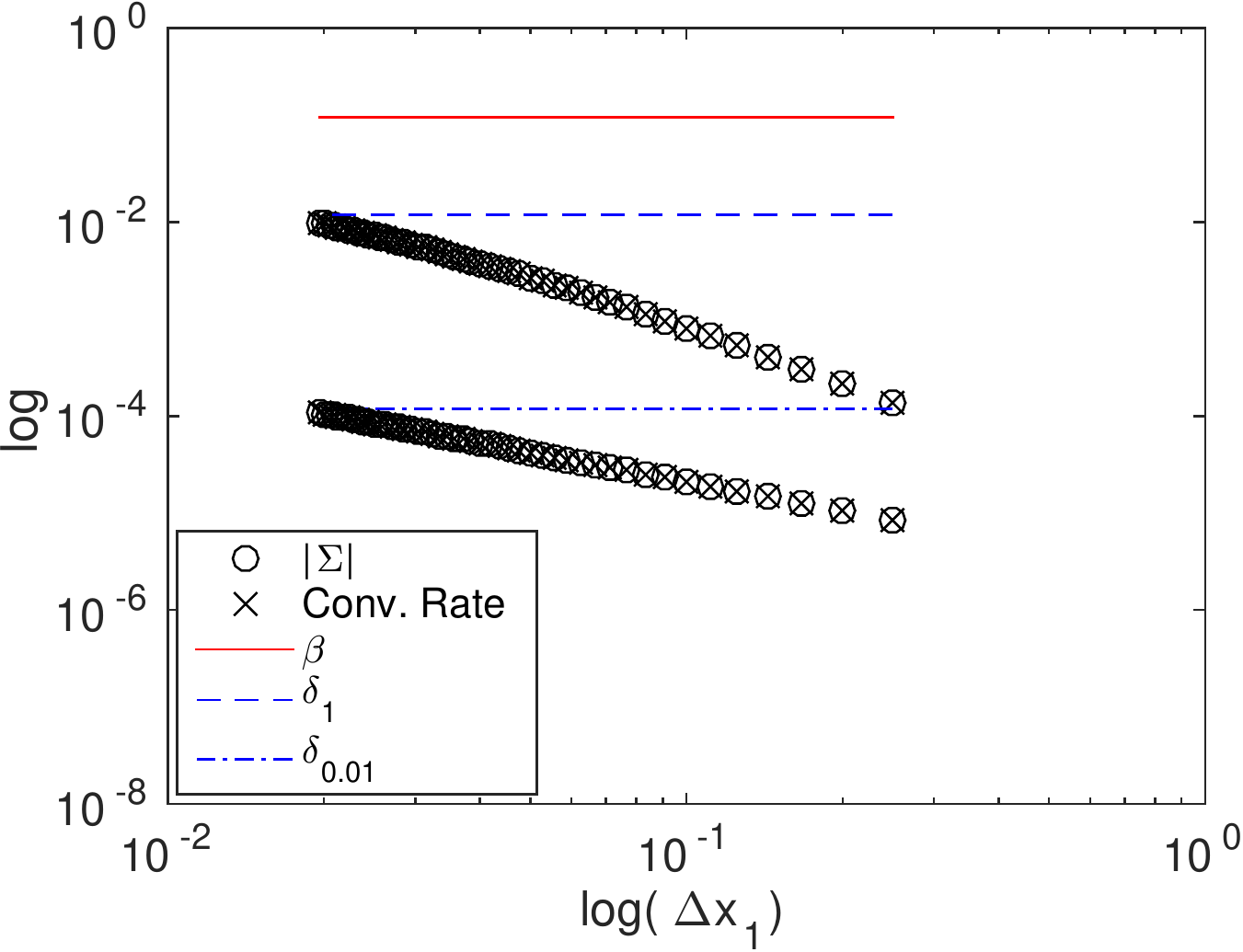}}
	\caption{Water-Steel thermal interaction with respect $\Delta t$ on the left and $\Delta x_1$ on the right in 1D.}
	\label{fig:water_steel}
\end{figure}

Figures \ref{fig:air_steel} and \ref{fig:water_steel} show the convergence rates for the interactions between air and steel and between water and steel, respectively. On the left we always have fixed $\Delta x_1$ and $r$, but variable $\Delta t$, whereas on the right we have fixed $\Delta t$ and $r$, but varying $\Delta x_1$. Each plot includes graphs for two different values of $r$. In figure \ref{fig:air_steel} we choose $r=1$ and $r=100$ to illustrate the effect of a neutral or a high aspect ratio. In figure \ref{fig:water_steel} we use $r=0.01$ and $r=1$ to illustrate how the rates are affected by a small or a neutral aspect ratio.  

Again, $|\Sigma|$ gives the exact convergence rates. Moreover, one observes that the rates in \ref{fig:air_steel}a and \ref{fig:water_steel}a tend to 0 as predicted in \eqref{limitdeltat1D} and in \ref{fig:air_steel}b and \ref{fig:water_steel}b to $\delta_r$ as predicted in \eqref{finallimitdeltax1D}. Furthermore, there is a roughly proportional relation between the convergence rate and the aspect ratio. For coupling with compressible flows, we typically have a high aspect ratio and therefore, the Dirichlet-Neumann iteration will be slowed down. Furthermore, this shows that it is very important to take the aspect ratio into account to make a reasonable prediction of the convergence rate at all. 

\subsection{2D FVM-FEM Results}

We now want to demonstrate that the 1D formula \eqref{formula1D} is a good estimator for the convergence rates in 2D. Thus, we now consider a 2D version of \eqref{model:eq} consisting of two coupled linear heat equations on two identical unit squares, e.g, $\Omega_1 = [-1,0] \times [0,1]$ and $\Omega_2 = [0,1] \times [0,1]$. We use a non equidistant cartesian grid with aspect ratio $r$ on $\Omega_1$ and an equidistant grid on $\Omega_2$. In order to use \eqref{formula1D} as an estimator we decided to take the equidistant mesh width on $\Omega_2$ as $\Delta x_2$ and the mesh width in $x$-direction on $\Omega_1$ as $\Delta x_1$. 
%For simplicity, we define $\Delta y := \Delta y_1 = \Delta y_2 = \Delta x_2 = 1/(N_2 +1)$ and $\Delta x := \Delta x_1 = 1/(N_1+1)$ with $N_1 \neq N_2$. 

As before, we present two real data examples described in table~\ref{tab:physicalcoeff} and~\ref{tab:asymptotics}, namely the thermal interaction between air at $273K$ with steel at $900K$ and air at $273K$ with water at $283K$. 

\begin{figure}[h!]
	\centering
	\subfigure[$\Delta t = 40/39, 2 \cdot 40/39,..., 39 \cdot 40/39$, $\Delta x_1 = 1/1100$ and $r = 100$ (top curves) or $r=1$ (bottom curves).]{\includegraphics[width=5.9cm]{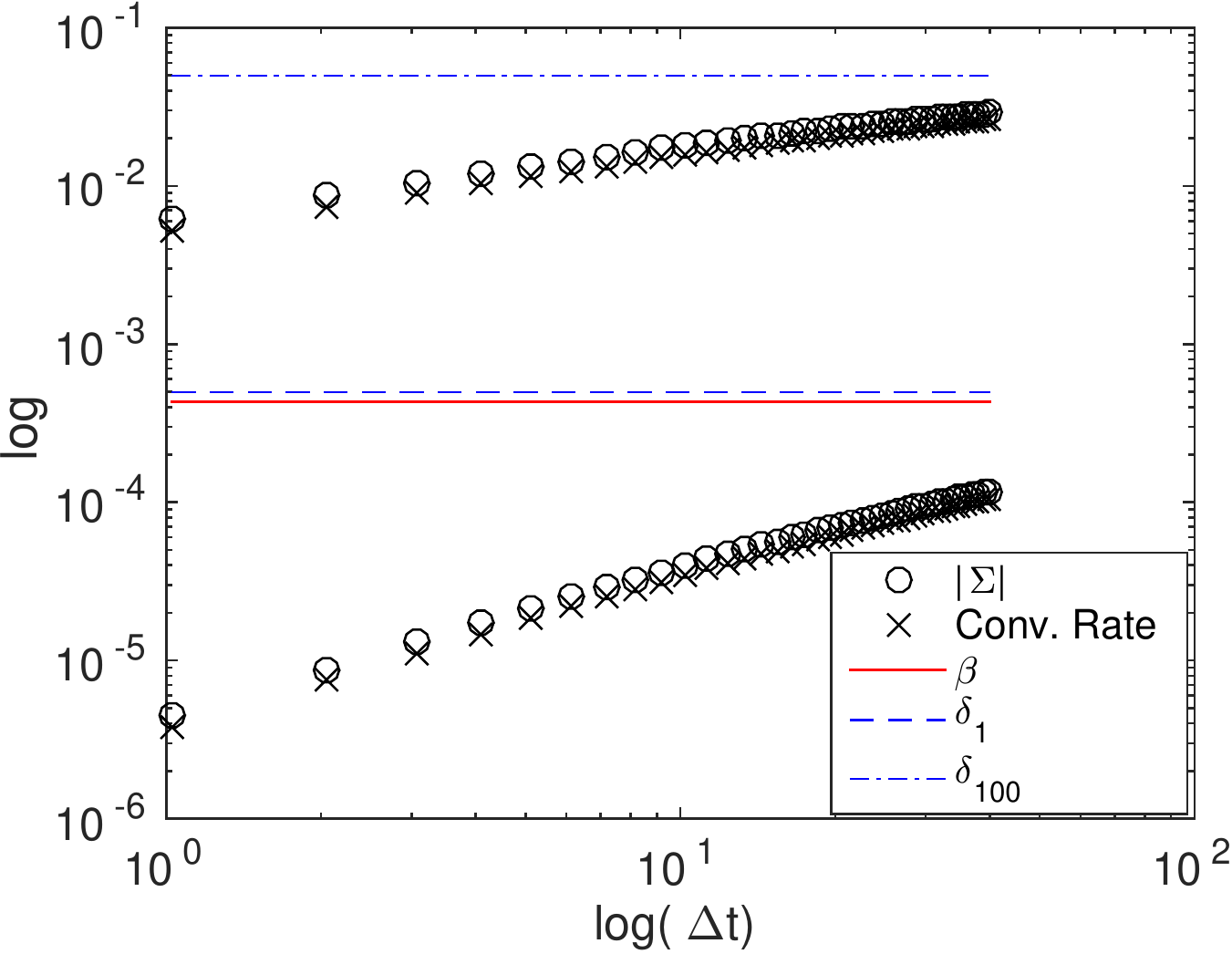}} \hfill
	\hspace{0.1cm}
	\subfigure[$\Delta x_1 = 1/3, 1/4,..., 1/50$, $\Delta t = 10$ and $r=100$ (top curves) or $r=1$ (bottom curves).]{\includegraphics[width=5.9cm]{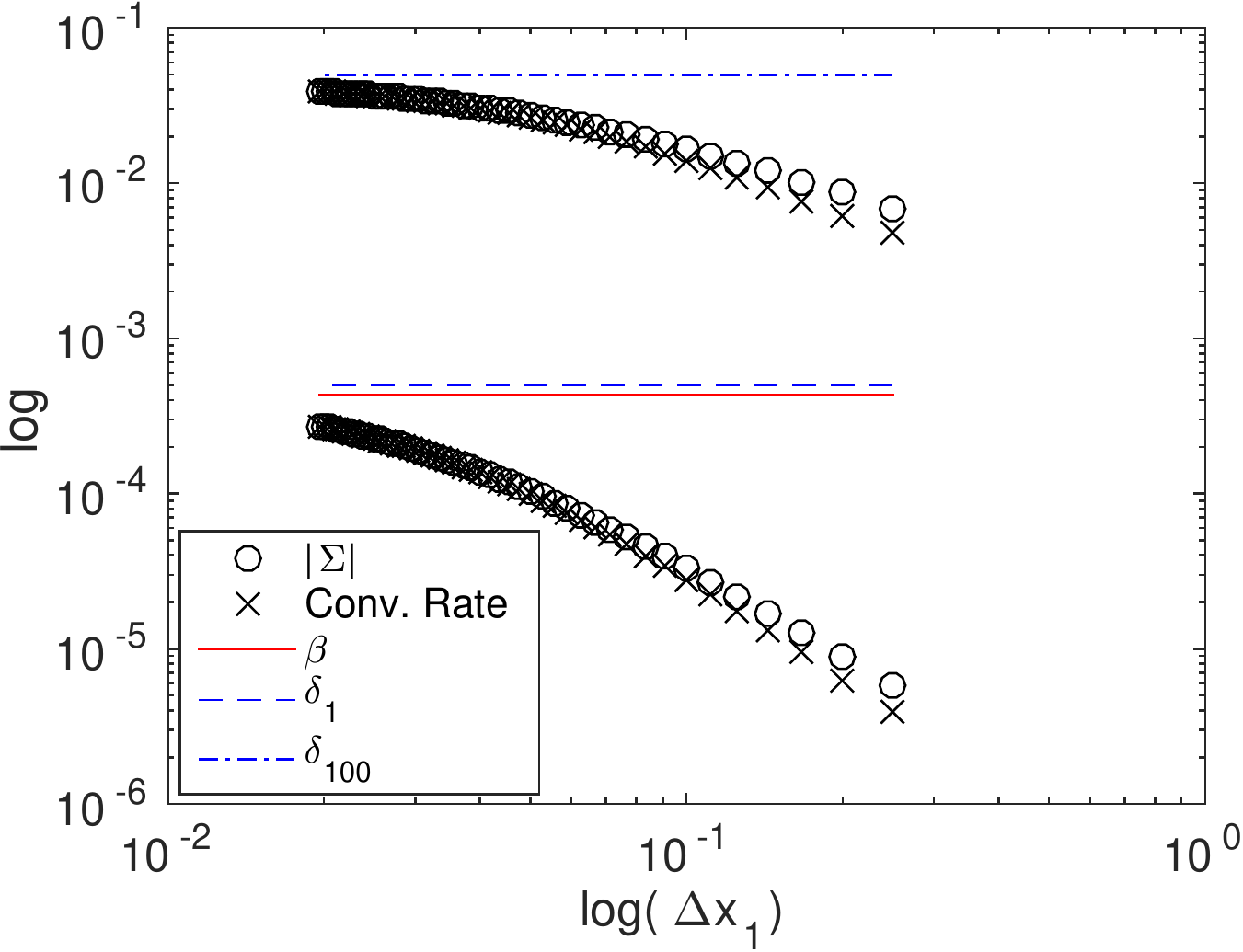}}
	\caption{2D Air-Steel thermal interaction. Observed and estimated convergence rates over $\Delta t$ (left) and $\Delta x_1$ (right).}
	\label{fig:air_steel2D}
\end{figure}

\begin{figure}[h!]
	\centering
	\subfigure[$\Delta t = 40/39, 2 \cdot 40/39,..., 39 \cdot 40/39$, $\Delta x_1 = 1/1100$ and $r = 1000$ (top curves) or $r=1$ (bottom curves).]{\includegraphics[width=5.9cm]{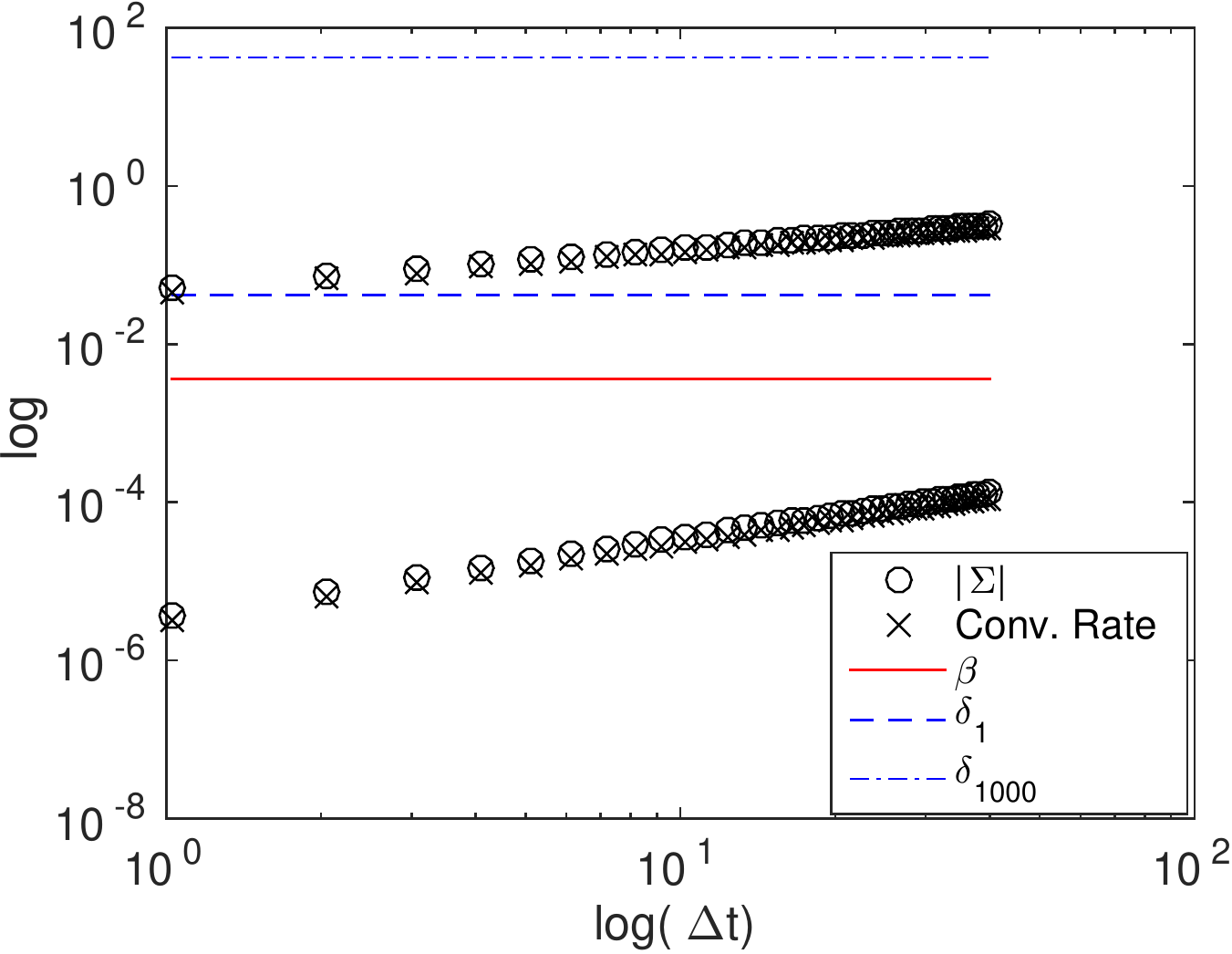}} \hfill
	\hspace{0.1cm}
	\subfigure[$\Delta x_1 = 1/3, 1/4,..., 1/35$, $\Delta t = 10$ and $r=1000$ (top curves) or $r=1$ (bottom curves).]{\includegraphics[width=5.9cm]{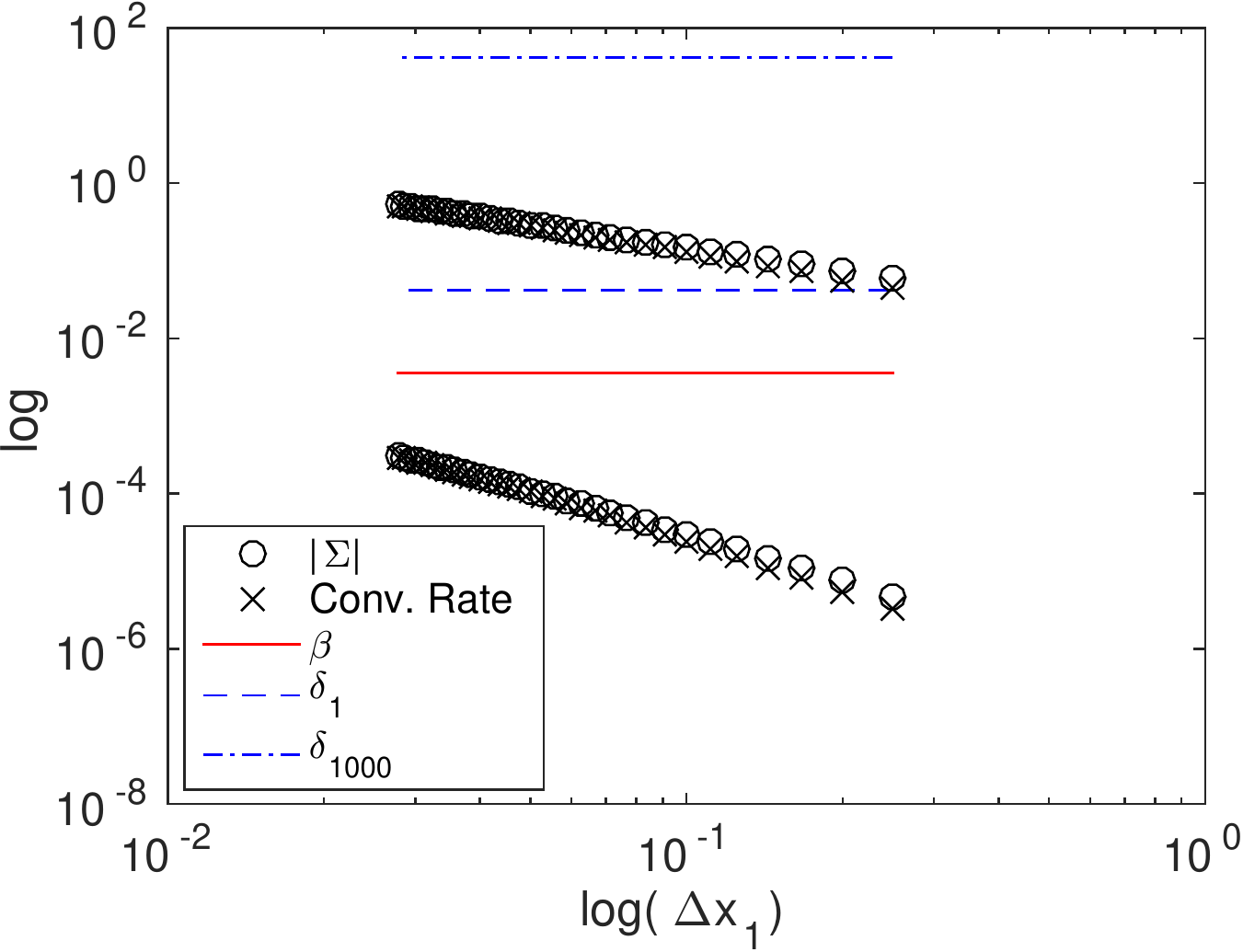}}
	\caption{2D Air-Water thermal interaction. Observed and estimated convergence rates over $\Delta t$ (left) and $\Delta x_1$ (right).}
	\label{fig:water_steel2D}
\end{figure}

Figures \ref{fig:air_steel2D} and \ref{fig:water_steel2D} show the convergence rates for the interactions between air and steel and between air and water in 2D respectively. On the left we always plot the rates for fixed $\Delta x_1$ and $r$ with variable $\Delta t$. On the right we plot the behaviour of the rates for fixed $\Delta t$ and $r$ and varying $\Delta x_1$. As before, notice that each plot includes two cases for two different $r$ values. In figure \ref{fig:air_steel2D} we choose $r=1$ and $r=100$ as in the 1D case (see figure \ref{fig:air_steel}) and in figure \ref{fig:water_steel2D} we use $r=1$ and $r=1000$ to illustrate the effect of a neutral or a high aspect ratio. One observes that the convergence rates predicted by the one-dimensional formula \eqref{formula1D} are almost exactly the ones observed in 2D. Thus, the 1D model problem case gives a very good estimator for the 2D model problem. 
%Moreover, one can claim when comparing figures \ref{fig:air_steel2D} and \ref{fig:air_steel} that the convergence behaviour of the Dirichlet-Neumann iteration in 1D and 2D under the same input data is almost equal. 

\subsection{Thermal FSI Test Cases}

\begin{figure}[h!]
	\centering
	\subfigure[Test case 1: Flow over a plate]{\includegraphics[width=6cm]{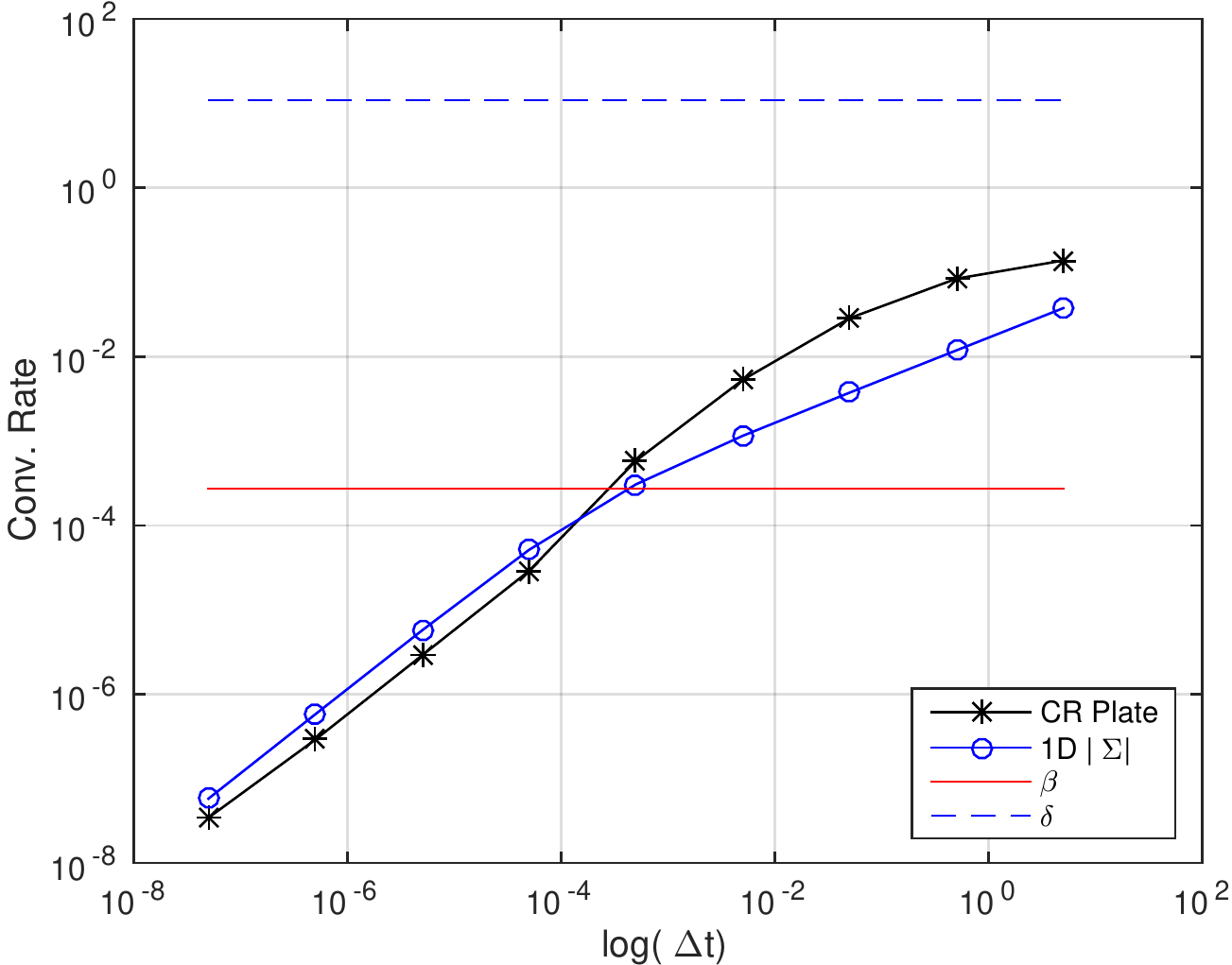}} \hfill
	\subfigure[Test case 2: Cooling of a flanged shaft]{\includegraphics[width=6cm]{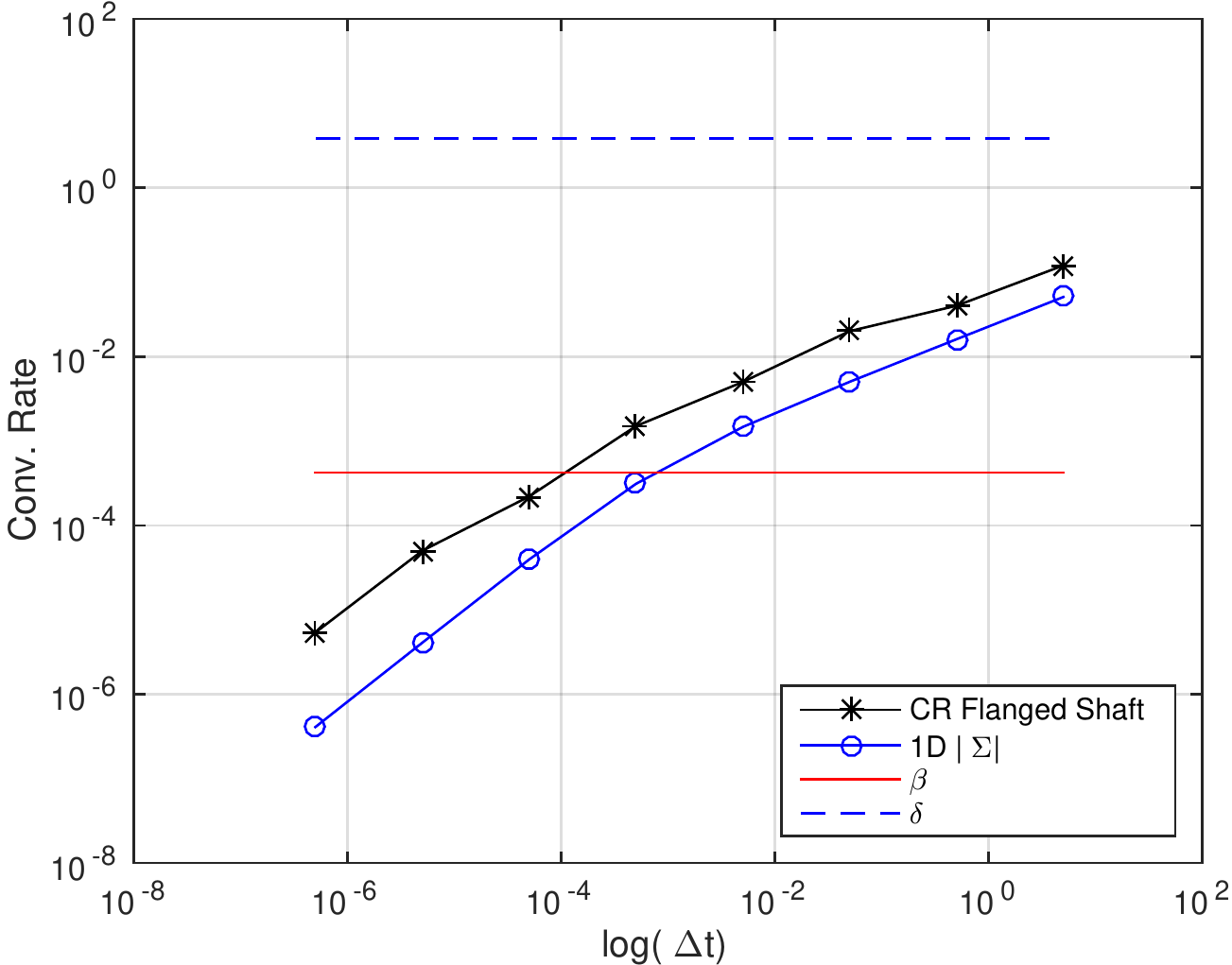}}
	\caption{Convergence behavior of the cooling systems with respect to $\Delta t$.}
	\label{fig:coolingsystemsestimator}
\end{figure}

Finally, we want to relate the results for the two nonlinear applications (the two cooling systems introduced in sections 3.1 and 3.2: the cooling of a flat plate and of a flanged shaft) to our analysis. Figure \ref{fig:coolingsystemsestimator}a shows the convergence behaviour for the flat plate and \ref{fig:coolingsystemsestimator}b for the flanged shaft. We plot the experimental convergence rates, the one-dimensional formula \eqref{formula1D}, the semidiscrete estimator \eqref{semidiscreteapprox} and the spatial limit $\delta_r$ specified in \eqref{finallimitdeltax1D}. 

In order to apply the 1D formula \eqref{formula1D} here, some assumptions need to be taken, since we partly have unstructured meshes and nonuniform temperatures. Thus, we assume air at $273K$ on the first subdomain with steel at $900K$ on the second subdomain for the cooling of a flat plate and air at $273K$ with steel at $1145K$ for the cooling of a flanged shaft. The density, heat capacity and heat conductivity of air and the density of steel are given in table~\ref{tab:physicalcoeff}. In addition, the heat conductivities and heat capacities of steel at $900K$ and $1145K$ are obtained from the nonlinear coefficient functions \eqref{heatconductivitysteel} and \eqref{heatcapacitysteel} by inserting $\Theta = 900K$ or $\Theta = 1145K$ respectively. This gives $\lambda = 39.82$ and $c_p = 1.3684e3$ for steel at $900K$ and $\lambda = 39.8$ and $c_p = 572.75$ for steel at $1145K$. 

Furthermore, for the cooling of a flat plate, we take $\Delta x_1 = 9.3736e-5$ which is the width of the fluid cells touching the interface in the $y$-direction and $\Delta x_2 = 1.6667$ which is width of the structure cells in both directions. Thus, we have an aspect ratio of $r=1.7780e4$. On the other hand, choosing $\Delta x_1$ and $\Delta x_2$ for the cooling of a flanged shaft is more difficult due to the unstructured grids. In order to get an upper bound for the aspect ratio $r$, we choose $\Delta x_1 = 1.6538e-4$ which is the minimum width of all the fluid cells touching the interface in direction normal to the wall and $\Delta x_2 = 1.1364$ which is the maximum width of all the structure cells touching the interface tangential to the wall. This gives $r=6.8713e3$. 

From figure \ref{fig:coolingsystemsestimator}a one observes with these choices, \eqref{formula1D} predicts the rates accurately for the cooling of a flat plate. Note that the semidiscrete estimator $\beta$ does not show any change with $\Delta t$. Remember that $\beta$ is almost always constant, except for a short dynamic transition between $(\lambda_1 /\lambda_2) \sqrt{D_2/D_1}$ and $\lambda_1/\lambda_2$ as shown in figure \ref{fig:semidiscretevsdeltat}. Here, we would have to choose a $\Delta t$ larger than $10E6$ to see the transition. 

Finally, in figure \ref{fig:coolingsystemsestimator}b one can see that \eqref{formula1D} predicts the convergence rates for the cooling of a flanged shaft to be only slightly smaller compared to the actual performance. This could be due to either the unstructured grids used or to the nonconstant temperature in the structure, which varies from room temperature to $1145K$. Again, $\beta$ is almost constant.

\section{Summary and Conclusions} 

We considered the Dirichlet-Neumann iteration for thermal FSI and studied the convergence rates. To this end, we considered the coupling of two heat equations on two identical domains. We assumed structured grids on both subdomains, but allowed for high aspect ratio grids in one domain. An exact formula for the convergence rates was derived for the 1D case. Furthermore, we determined the limits of the convergence rates when approaching the continuous case either in space ($r\lambda_1 / \lambda_2 $) or time ($0$). This was confirmed by numerical results, where we also demonstrated that the 1D case gives excellent estimates for the 2D case. In addition, numerical experiments show that the linear analysis is relevant for nonlinear thermal FSI problems.  

All in all, strong jumps in the coefficients of the coupled PDEs will imply fast convergence. In the domain decomposition context, the coupling will be slow because the material coefficients are continuous over all the subdomains, i.e, $\lambda_1 = \lambda_2$, and therefore $\delta_1 \sim 1$. For coupling of structures and compressible flows, the aspect ratio in the structure has to be taken into account, since the convergence rate is proportional to it. For the nonlinear cooling problems considered here, the convergence rate was still around 0.1 for large $\Delta t$. When encountering divergence anyhow, this can be solved by reducing the time step. Note that in a time adaptive setting, it is standard to allow for a feedback loop between the nonlinear solver and the time stepper. 

\bibliography{birkenmonge}

\end{document}